\documentclass[11pt,a4paper]{amsart}

\usepackage{amsmath}
\usepackage{amssymb}
\usepackage{enumerate}
\usepackage[final]{graphicx}
\usepackage[papersize={165mm,235mm},top=26mm,left=16mm,text={128mm,192mm},%
footskip=4pt,headsep=6mm,marginparsep=0pt]{geometry}
\usepackage{url}
\usepackage{ulem, color}
\newtheorem{theorem}{Theorem}[section]
\newtheorem{corollary}[theorem]{Corollary}
\newtheorem{lemma}[theorem]{Lemma}
\newtheorem{proposition}[theorem]{Proposition}
\theoremstyle{definition}
\newtheorem{definition}[theorem]{Definition}
\newtheorem{example}[theorem]{Example}
\newtheorem{remark}[theorem]{Remark}

\numberwithin{equation}{section}

\title[$UC$ and $UC^{*}$ properties and best proximity points]{On the $UC$ and $UC^{*}$ properties and the existence of best proximity points in metric spaces}

\author[V. Zhelinski]{Vasil Zhelinski}

\address[V. Zhelinski]{University of Plovdiv Paisii Hilendarski,
	24 Tsar Assen str, 4000 Plovdiv, Bulgaria}
\email{{\tt vasil\_zhelinski@uni-plovdiv.bg}}

\author[B. Zlatanov]{Boyan Zlatanov}
\address[B. Zlatanov]{University of Plovdiv Paisii Hilendarski,
	24 Tsar Assen str, 4000 Plovdiv, Bulgaria}
\email{\tt bobbyz@uni-plovdiv.bg}

\begin{document}

\begin{abstract}
We investigate the connections between $UC$ and $UC^{*}$ properties for ordered pairs of subsets $(A,B)$ in metric spaces, which are involved in the study of existence and uniqueness of best proximity points. We show that the $UC^{*}$ property is included into the $UC$ property. We introduce some new notions: bounded $UC$ ($BUC$) property and uniformly convex set about a function $\phi$.
We prove that these new notions are generalizations of the $UC$ property and that both of them are sufficient for to ensure existence and uniqueness of best proximity points. We show that these two new notions are different from a uniform convexity and even from a strict convexity. If we consider the underlying space to be a Banach space we find a sufficient condition which ensures that from the $UC$ property it follows   the uniform convexity of the underlying Banach space. We illustrate the new notions with examples. We present an example of a cyclic contraction $T$ in a space, which is not even strictly convex and the ordered pair $(A,B)$ has not the $UC$ property, but has the $BUC$ property and thus there is a unique best proximity point of $T$ in $A$.\\




\end{abstract}

\maketitle


\section{Introduction}

A fundamental result in the fixed point theory is the Banach contraction principle in Banach spaces or in complete metric spaces. Fixed point theory is an important tool for solving equations $Tx=x$ for mapping T defined on subsets of metric or normed spaces. 

One kind of a generalization of the Banach contraction principle is the notion of cyclic maps, $T:A\to B$ and $T:B\to A$ \cite{Kirk-Srinivasan-Veeramani}. Because a non-self mapping $T:A\to B$ does not necessarily have a fixed point, one often attempts to find an element $x$ which is in some sense ''closest`` to $Tx$. 
Thus we can alter the fixed point problem into the optimization problem $\min\{\|x-Tx\|:x\in A\cup B\}$. Best proximity point theorems, introduced in \cite{Eldred-Veeramani}, are relevant in this perspective. 
A sufficient condition for the existence and uniqueness of the best proximity points in uniformly convex Banach spaces is given in \cite{Eldred-Veeramani}. The uniform convexity of the underlying space ensures good geometric properties of the space { and is a key property in getting the results of existence and uniqueness of best proximity points}.

Naturally \cite{Sintunavarat-Kumam,Suzuki-Kikkawa-Vetro}, {in} solving the optimization problem $\min\{\|x-Tx\|:x\in A\cup B\}$ for a cyclic map $T:A\to B$ and $T:B\to A$, where $A,B$ are subsets of { either} a Banach space $(X,\|\cdot\|)$ { or a metric space $(X,\rho)$},
{ it} may {be} need{ed} only some specific properties of the domain of $T$, i.e. $A\cup B$, instead of the uniform convexity of the underlying Banach space { $(X,\|\cdot\|)$}.

This idea to search for good properties of the {ordered pair $(A,B)$ of} sets, which defines the domain of the cyclic map, have been firstly initiated in \cite{Suzuki-Kikkawa-Vetro}, where the authors have investigated existence and uniqueness of best proximity points in metric space over ordered pairs $(A,B)$ of sets in complete metric spaces. The authors have introduced the notion for an ordered pair $(A,B)$ of sets to satisfy the property $UC$. Some relations of the property $UC$ for sets in a Banach space and the properties of uniform convexity, uniform convexity in every direction and relatively compact sets were presented in \cite{Suzuki-Kikkawa-Vetro}.

Later on a new notion for an ordered pair $(A,B)$ of sets to satisfy the property $UC^{*}$ was introduced, in order to investigate existence and uniqueness of coupled best proximity points in complete metric spaces, rather than uniformly convex Banach spaces \cite{Sintunavarat-Kumam}.

{Deep results about fixed points and the geometry of the underlying space can be found in \cite{Berinde,Chidume,Khamsi-Kozlowski}.}

Some results about applications of coupled best proximity points for solving of symmetric \cite{Ilchev-Zlatanov} and for solving of non-symmetric \cite{Zlatanov} systems of equations have been presented. It is interesting to mention that the presented technique in \cite{Ilchev-Zlatanov,Zlatanov} enables to find exact solutions in cases, where the classical {fixed point} methods can find only approximations. Best proximity point results have been used in searching of market equilibrium in duopoly markets, where the cyclic maps have been replaced by semi-cyclic maps \cite{Ajeti-Ilchev-Zlatanov,Dzhabarova-Kabaivanov-Ruseva-Zlatanov}. The natural underlying space in the market equilibrium theory is close to non convex spaces, rather than convex spaces as pointed in \cite{Ajeti-Ilchev-Zlatanov}, where results about coupled best proximity points have been obtained in reflexive Banach spaces.

Thereafter it seem interesting to search for some conditions, different from the uniform convexity of the underlying space, that will ensure {some of} the properties, involved in the definitions of $UC$ or/and $UC^{*}$ {and will lead to positive conclusions on the existence and uniqueness of best proximity points}.

\section{Preliminaries}

In what follows we will use the notations: $\mathbb{N}$ for the set of natural numbers, $\mathbb{R}$ for the set of real numbers, $S_X$ and $B_X$ for the unit sphere and the unit ball, respectively, where $(X,\|\cdot\|)$ is a Banach space, $B(x_0,r)=\{x\in X:\|x-x_0\|<r$\} and $B[x_0,r]=\{x\in X:\|x-x_0\|\leq r$\} for the open and close balls with a center $x_0$ and radius $r$, respectively. Let $(X,\rho)$ be a metric space and $A,B\subset X$. We will denote by ${\rm dist}(A,B)=\inf\{\rho(a,b):A\in A, b\in B\}$ the distance between the sets $A$ and $B$. {Whenever} the underlying space is a Banach space $(X,\|\cdot\|)$ we will consider the metric to be the one generated by the norm, i.e. $\rho(x,y)=\|x-y\|$.

\begin{definition}\label{Definition:1}(\cite{Clarkson,FHHMZ})
	Let $(X,\|\cdot\|)$ be a Banach space. For every $\varepsilon\in (0,2]$ we define the modulus of convexity of $\|\cdot\|$ by
	$$
	\delta_{\|\cdot\|}(\varepsilon)=\inf\left\{1-\left\|\frac{x+y}{2}\right\|:x,y\in B_X, \|x-y\|\geq\varepsilon\right\}.
	$$
	The norm is called uniformly convex if $\delta_X(\varepsilon)>0$ for all $\varepsilon\in (0,2]$. The space $(X,\|\cdot\|)$ is then called uniformly convex {Banach} space.
\end{definition}

\begin{definition}\label{Definition:2}(\cite{Kirk-Srinivasan-Veeramani})
	Let $A$ and $B$ be two sets. A map $T:A\cup B\to A\cup B$ is called a cyclic map if it satisfies $T:A\to B$ and $T:B\to A$. 
\end{definition}

\begin{definition}\label{Definition:3}(\cite{Eldred-Veeramani})
	Let $(X,\rho)$ be a metric space, $A$ and $B$ be subsets of $X$ and $T:A\cup B\to A\cup B$ be a cyclic map. 
	We say that the point $x\in A$ is a best proximity point of $T$ in $A$, if $\rho(x,Tx)={\rm dist}(A,B)$. 
\end{definition}

\begin{definition}\label{Definition:4}(\cite{Eldred-Veeramani})
	Let $(X,\rho)$ be a metric space, $A$ and $B$ be subsets of $X$. We say that the map $T:A\cup B\to A\cup B$ is a cyclic contraction map, if it is a cyclic map and satisfies the inequality
	$$
	\rho(Tx,Ty)\leq \alpha\rho(x,y)+(1-\alpha){\rm dist}(A,B)
	$$
	for some $\alpha \in (0,1)$ and every $x\in A$, $y\in B$.
\end{definition}

\begin{theorem}\label{Theorem:5}(\cite{Eldred-Veeramani})
	Let $A$ and $B$ be nonempty closed and convex subsets of a uniformly convex
	Banach space $(X,\|\cdot\|)$. Suppose $T:A\cup B\to A\cup B$ be a cyclic contraction map, then there exists a
	unique best proximity point $x$ of $T$ in $A$.
\end{theorem}

It is also proven in \cite{Eldred-Veeramani}, that for any initial guess $x_0\in A${,} the iterated sequence $x_n=T^nx_0$, for $n\in \mathbb{N}${,} splits into two sequences, such that
$\left\{x_{2n}\right\}_{n=1}^\infty$ converges to the best proximity point $x$ of $T$ in $A$ and $\left\{x_{2n-1}\right\}_{n=1}^\infty$ converges to the best proximity point $Tx$ of $T$ in $B$. The apriori and the aposteriori error estimates have been found \cite{Zlatanov-0} of the iterated sequences $\{x_{2n}\}_{n=1}^\infty$ and $\{x_{2n-1}\}_{n=1}^\infty$.

The next lemmas are crucial in getting of the results about best proximity points in uniformly convex Banach spaces.

\begin{lemma}\label{Lemma:6}(\cite{Eldred-Veeramani})
	Let $(X,\|\cdot\|)$ be a uniformly convex Banach space. Let $A$ and $B$ be nonempty and closed subsets of $X$. Let $A$ be convex. Let $\{ x_n\}_{n=1}^\infty$ and $\{ z_n\}_{n=1}^\infty$ be sequences in $A$ and $\{ y_n \}_{n=1}^\infty$ be a sequence in $B$ such that $\displaystyle\lim_{{n\to\infty}}\|x_n - y_n\|=\lim_{{n\to\infty}}\|z_n - y_n\|={\rm dist}(A,B)$. Then $\displaystyle\lim_{{n\to\infty}} \|x_n - z_n\| = 0$.
\end{lemma}

\begin{lemma}\label{Lemma:7}(\cite{Eldred-Veeramani})
	Let $X$ be a uniformly convex Banach space. Let $A$ and $B$ be non empty and closed subsets of $X$. Let $A$ be convex. Let $\{ x_n\}_{n=1}^\infty$ and $\{ z_n\}_{n=1}^\infty$ be sequences in $A$ and $\{ y_n \}_{n=1}^\infty$ be a sequence in $B$ such that:
	\begin{itemize}\setlength{\itemsep}{0em}
		\item $\lim_{n\to\infty}\|z_n- y_n\|={\rm dist}(A,B)$
		\item for every $\varepsilon>0$ there is $N\in\mathbb{N}$ so that the inequality
		$\|x_m-y_n\|\leq {\rm dist}(A,B)+\varepsilon$ holds for all $m>n\geq N$.
	\end{itemize}
	Then for every $\varepsilon>0$ there exists $N_0\in\mathbb{N}$ so that for all $m>n\geq N_0$ there holds the inequality 
	$\|x_m-z_n\|\leq \varepsilon$.
\end{lemma}

\begin{definition}\label{Definition:8}(\cite{FHHMZ})
	Let $(X,\|\cdot\|)$ be a Banach space. $X$ is called a strictly convex Banach space if $\|x+y\|<2$ for all $x,y\in S_X$, such that $x\not=y$.
\end{definition}

The next lemma is actually proven in Proposition 5 in \cite{Suzuki-Kikkawa-Vetro}, without stating it as a particular proposition.

\begin{lemma}\label{Lemma:9}(\cite{Suzuki-Kikkawa-Vetro}, Proposition 5)
	Let $A$, $B$ be closed subsets of a strictly convex normed space $(X,\|\cdot\|)$, such that ${\rm dist}(A,B)>0$ and let $A$ be convex. If
	$x, z\in A$ and $y\in B$ be such that $\|x-y\|=\|z-y\|={\rm dist}(A,B)$, then $x=z$.
\end{lemma}

As pointed in \cite{Eldred-Veeramani}, if the sets $A$ and $B$ satisfy some additional properties, $T:A\cup B\to A\cup B$ be a cyclic contraction map and either $A$ or $B$ is boundedly compact, then there exists a best proximity point $x$ of $T$ in $A$.

The authors of \cite{Suzuki-Kikkawa-Vetro} have found some properties of the sets $A$ and $B$ that define the domain of the cyclic contraction map $T$, which ensure the existence and uniqueness of the best proximity points, without assuming for the underlying space to be a uniform convex Banach space.

\begin{definition}\label{Definition:10}(\cite{Suzuki-Kikkawa-Vetro}) 
	Let $A$ and $B$ be nonempty subsets of a metric space $(X,\rho)$. We say that the ordered pair $(A,B)$ 
	satisfies the property $UC$ if for every sequences $\left\{x_n\right\}_{n=1}^\infty, \left\{z_n\right\}_{n=1}^\infty\subset A$ and
	$\left\{y_n\right\}_{n=1}^\infty\subset B$, such that there holds $\displaystyle\lim_{n\to\infty}\rho\left(x_n,y_n\right)=\lim_{n\to\infty}\rho\left(z_n,y_n\right)={\rm dist}(A,B)$, then there holds
	$\displaystyle\lim_{n\to\infty}\rho\left(x_n,z_n\right)=0$.
\end{definition}

It is easy to observe that the $UC$ property replaces Lemma \ref{Lemma:6} and that the assumption in Lemma \ref{Lemma:6} for the sets $A$ and $B$ to be closed ones and $A$ to be convex is not necessary, when we replace the uniform convexity of the underlying Banach space with the $UC$ property. We would like to points out that it may happen the ordered pair $(A,B)$ to satisfy $UC$, but the ordered pair $(B,A)$ not to satisfy it.

Some property of the $UC$ ordered pairs $(A,B)$ of subsets are presented in \cite{Suzuki-Kikkawa-Vetro}.

\begin{proposition}\label{Proposition:11}(\cite{Suzuki-Kikkawa-Vetro}) 
	Let $A$ and $B$ be nonempty subsets of a uniformly convex Banach space $(X,\|\cdot\|)$. 
	If $A$ is convex, then the ordered pair $(A,B)$ has the property $UC$.
\end{proposition}

\begin{proposition}\label{Proposition:12}(\cite{Suzuki-Kikkawa-Vetro}) 
	Let $A$ and $B$ be nonempty subsets of a metric space $(X,\rho)$, such that ${\rm dist}(A,B)=0$.
	Then the ordered pair $(A,B)$ satisfies the property UC.
\end{proposition}

\begin{proposition}\label{Proposition:13}(\cite{Suzuki-Kikkawa-Vetro}) 
	Let $A, A^\prime, B, B^\prime$ be nonempty subsets of a metric space 
	$(X,\rho)$, such that $A\subseteq A^\prime$, $B\subseteq B^\prime$  and ${\rm dist}(A,B)={\rm dist}(A^\prime,B^\prime)$. If the ordered pair $(A^\prime, B^\prime$) satisfies the property UC, then the ordered pair $(A,B)$ satisfies the property UC too.
\end{proposition}

If in the next result the underlying metric space $(X,\rho)$ is replaced by a uniformly convex Banach space $(X,\|\cdot\|)$ and $A$ and $B$ be nonempty closed and convex subsets then it generalizes Theorem \ref{Theorem:5}.

\begin{theorem}\label{Theorem:14}(\cite{Suzuki-Kikkawa-Vetro}) 
	Let $A$ and $B$ be nonempty closed subsets of a complete metric space $(X,\rho)$,
	such that the ordered pairs $(A, B)$ satisfy the property $UC$. 
	Let $T:A\cup B\to A\cup B$ be a cyclic map and there exists $k\in [0,1)$, so that the inequality
	$$
	\rho(Tx,Ty)\leq k\max\{\rho(x,y),\rho(x,Tx),\rho(y,Ty)\}+(1-k){\rm dist}(A,B)
	$$
	holds for all $x\in A$ and $y\in B$. Then there is a unique best proximity point $x$ of $T$ in $A$, the sequence of successive iterations $\{T^{2n}x_0\}_{n=1}^\infty$ converges to $x$ for any initial guess $x_0\in A$. There is at least one best proximity point $y\in B$ of $T$ in $B$.  Moreover best proximity point $y\in B$ of $T$ in $B$ is unique, provided that the ordered pair $(B,A)$ has the $UC$ property.
\end{theorem}

Examples show that it is possible to have a Banach space $(X,\|\cdot\|)$ that is not even strictly convex, but there are sets $A$ and $B$, so that the ordered pair $(A,B)$ satisfies the property $UC$. 

\begin{figure}[h!]
	\centering
	\includegraphics [width=0.80\textwidth] {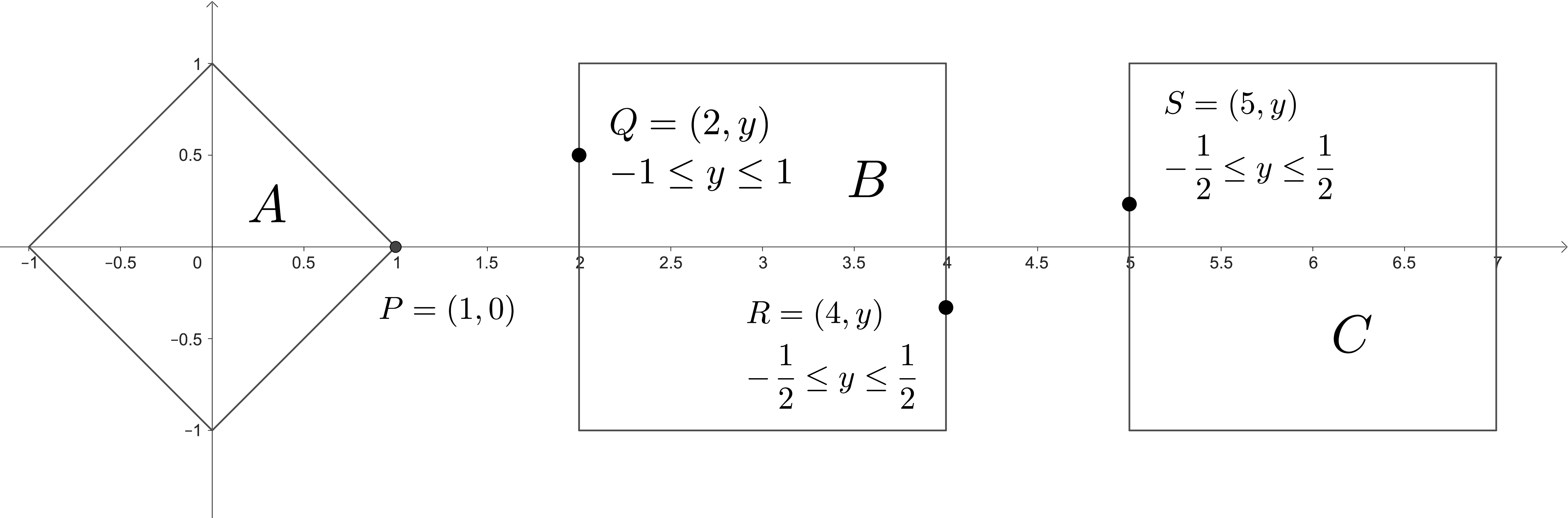}
	\caption{Example \ref{Example:15}}  \label{Fig01}
\end{figure}

\begin{example}\label{Example:15}
	Let us consider the spaces $X_\infty=(\mathbb{R}^2,\|\cdot\|_\infty)$, $X_1=(\mathbb{R}^2,\|\cdot\|_1)$ and the sets 
	$A,B\subset X_1$ and $B,C\subset X_\infty$ (Figure \ref{Fig01}). 
\end{example}

It easy to see that ${\rm dist}(A,B)=\inf\{\|x-y\|_1:x\in A,y\in B\}=1$.

If {the} sequences $\{a_n\}_{n=1}^\infty\subset A$ and $\{b_n\}_{n=1}^\infty\subset B$
satisfy $\lim_{n\to\infty}\|a_n-b_n\|={\rm dist}(A,B)$, then $\lim_{n\to\infty}a_n=(1,0)$ and $\lim_{n\to\infty}b_n=(2,0)$.
Therefore if there hold $\lim_{n\to\infty}\|x_n-y_n\|={\rm dist}(A,B)$ and $\lim_{n\to\infty}\|z_n-y_n\|={\rm dist}(A,B)$ for
$\{x_n\}_{n=1}^\infty,\{z_n\}_{n=1}^\infty\subset A$ and $\{y_n\}_{n=1}^\infty\subset B$, then there holds 
$\lim_{n\to\infty}\|x_n-z_n\|=0$. Thus the ordered pair of sets $(A,B)$ has the $UC$ property.

There holds ${\rm dist}(B,C)=\inf\{\|x-y\|_\infty:x\in B,y\in C\}=1$. 
Let us consider the sequences $\{c_n\}_{n=1}^\infty\subset C$ and $\{b_n\}_{n=1}^\infty\subset B$
satisfy $\lim_{n\to\infty}\|b_n-c_n\|={\rm dist}(B,C)$, where $c_n=(x_n,y_n)$ and $b_n=(u_n,v_n)$, be such that $\lim_{n\to\infty}x_n=5$, $\lim_{n\to\infty}u_n=4$ and $y_n,v_n\in [-1/2,1/2]$ can be arbitrary.
Therefore there exist sequences $\{x_n\}_{n=1}^\infty,\{z_n\}_{n=1}^\infty\subset B$ and $\{y_n\}_{n=1}^\infty\subset C$, satisfying $\lim_{n\to\infty}\|x_n-y_n\|={\rm dist}(B,C)$ and $\lim_{n\to\infty}\|z_n-y_n\|={\rm dist}(B,C)$, so that the
$\lim_{n\to\infty}\|x_n-z_n\|$ does not exist. Consequently the ordered pair of sets $(B,C)$ has not the $UC$ property.

If we consider the set $\mathbb{R}^2$ endowed with a{color{red}\sout{n}} uniformly convex norm (for example the Hilbert norm $\|\cdot\|_2$, then 
the ordered pair of sets $(B,C)$ has the $UC$ property according to Lemma \ref{Lemma:6}.

\begin{figure}[h!]
	\centering
	\includegraphics [width=0.80\textwidth] {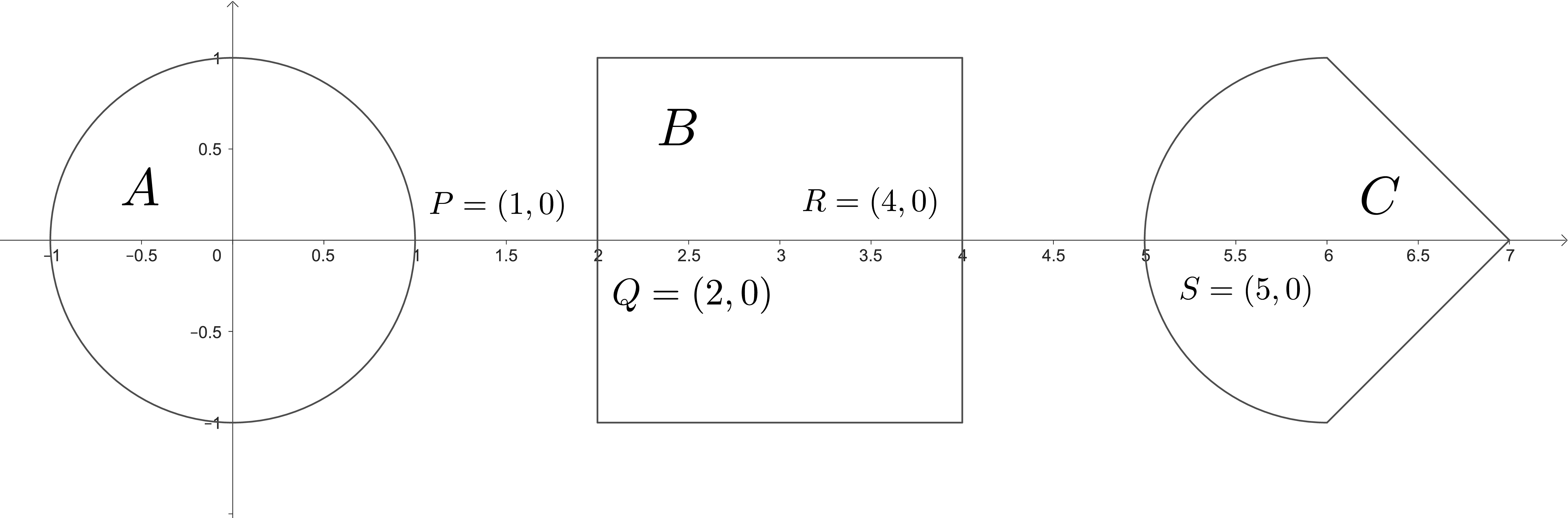}
	\caption{Example \ref{Example:16}}  \label{Fig02}
\end{figure}

\begin{example}\label{Example:16}
	Let us consider the space $X_\infty=(\mathbb{R}^2,\|\cdot\|_\infty)$ and the sets $A$, $B$ and $C$ (Figure \ref{Fig02}). 
\end{example}

It easy to see that ${\rm dist}(A,B)={\rm dist}(B,C)=1$. 

If {the} sequences $\{a_n\}_{n=1}^\infty\subset A$ and $\{b_n\}_{n=1}^\infty\subset B$
satisfy $\lim_{n\to\infty}\|a_n-b_n\|={\rm dist}(A,B)$, then $\lim_{n\to\infty}a_n=(1,0)$ and $\lim_{n\to\infty}b_n=(2,0)$.
Therefore if there hold $\lim_{n\to\infty}\|x_n-y_n\|={\rm dist}(A,B)$, $\lim_{n\to\infty}\|z_n-y_n\|={\rm dist}(A,B)$ for the sequences
$\{x_n\}_{n=1}^\infty,\{z_n\}_{n=1}^\infty\subset A$ and $\{y_n\}_{n=1}^\infty\subset B$, then there holds 
$\lim_{n\to\infty}\|x_n-z_n\|=0$. Thus the ordered pair of sets $(A,B)$ has the $UC$ property.

By similar arguments we get that the ordered pair of sets $(B,C)$ has the $UC$ property, never mind the geometry of the unit ball in $\mathbb{R}^2$.

By the observations, in the two examples, we see that the property $UC$ depend on three conditions: the geometry of the unit ball of the underlying Banach space, the geometry properties of the sets $A$ and $B$ and the positioning of the sets into the space.

\begin{definition}\label{Definition:17}(\cite{DJS})
	Let $(X,\|\cdot\|)$ be a Banach space. For every $\varepsilon\in (0,2]$ and every $z\in X\backslash\{0\}$ we define $\delta_{\|\cdot\|}(z,\varepsilon)$ by
	$$
	\delta_{\|\cdot\|}(z,\varepsilon)=\inf\left\{1-\textstyle\left\|\frac{x+y}{2}\right\|:x,y\in S_X, \|x-y\|\geq\varepsilon, x-y=\lambda z\ \mbox{for some}\ \lambda\in\mathbb{R}\right\}.
	$$
	We call $\delta_{\|\cdot\|}(z,\varepsilon)$ the modulus of convexity in the direction $z\in X\backslash\{0\}$.
	
	The norm is called uniformly convex in every direction ($UCED$) if $\delta_X(z,\varepsilon)>0$ for all $z\in X\backslash\{0\}$ and $\varepsilon\in (0,2]$. The space $(X,\|\cdot\|)$ is then called uniformly convex in every direction Banach space.
\end{definition}

\begin{proposition}\label{Proposition:18}(\cite{Suzuki-Kikkawa-Vetro}) 
	Let $A$ and $B$ be nonempty subsets of a $UCED$ Banach space $(X,\|\cdot\|)$. If $A$ is convex and relatively compact then the ordered pair $(A,B)$ has the property $UC$.
\end{proposition}

A notion that replaces Lemma \ref{Lemma:7} have been introduced in \cite{Sintunavarat-Kumam} in order to get existence and uniqueness results about coupled best proximity points in metric spaces.

\begin{definition}\label{Definition:19}(\cite{Sintunavarat-Kumam})
	Let $A$ and $B$ be nonempty subsets of the metric space $(X,\rho)$. We say that the ordered pair $(A,B)$ satisfies the property $UC^{*}$ if $(A,B)$ satisfies the property $UC$ and for every sequences $\{x_n\}_{n=1}^\infty, \{z_n\}_{n=1}^\infty\subset A$, $\{y_n\}_{n=1}^\infty\subset B$, so that
	\begin{enumerate}[\ref{Definition:19}.1)]
		\item\label{enum-1}
		$$
		\displaystyle\lim_{{n\to\infty}}\rho(z_n,y_n)={\rm dist}(A,B)
		$$
		\item\label{enum-2} for every $\varepsilon>0$ there exists $N\in\mathbb{N}$, so that the inequality  
		$$\rho(x_m,y_n)\leq {\rm dist}(A,B)+\varepsilon$$
		holds for all $m>n\geq N$.
	\end{enumerate}
	Then for any $\varepsilon>0$ there is $N_1\in\mathbb{N}$ so that the inequality  
	$\rho(x_m,z_n)\leq \varepsilon$ holds for all $m>n\geq N_1$.
\end{definition}

For easier presentation of the results we will alter Definition \ref{Definition:19} slightly.
	
	{\bf Definition: 19a}
	{\it Let $A$ and $B$ be nonempty subsets of the metric space $(X,\rho)$. We say that the ordered pair $(A,B)$ satisfies the property restricted $UC^{*}$ ($RUC^{*}$) if for every sequences $\{x_n\}_{n=1}^\infty, \{z_n\}_{n=1}^\infty\subset A$, $\{y_n\}_{n=1}^\infty\subset B$, so that
		\begin{itemize}
			\item
			$\displaystyle\lim_{{n\to\infty}}\rho(z_n,y_n)={\rm dist}(A,B)$
			\item for every $\varepsilon>0$ there exists $N\in\mathbb{N}$, so that the inequality  
			$\rho(x_m,y_n)\leq {\rm dist}(A,B)+\varepsilon$
			holds for all $m>n\geq N$.
		\end{itemize}
		Then for every $\varepsilon>0$ there is $N_1\in\mathbb{N}$ so that the inequality  
		$\rho(x_m,z_n)\leq \varepsilon$ holds for all $m>n\geq N_1$.}
	
	We have removed the assumption that the ordered pair $(A,B)$ satisfies the $UC$ property.

\begin{proposition}\label{Proposition:20}(\cite{Sintunavarat-Kumam})
	Every ordered pair $(A,B)$ of nonempty subsets of a metric space $(X,\rho)$, so that ${\rm dist}(A, B)=0$,
	satisfies the property $UC^*$.
\end{proposition}

\begin{proposition}\label{Proposition:21}(\cite{Sintunavarat-Kumam})
	Every ordered pair $(A,B)$ of nonempty subsets of a uniformly convex Banach space $(X,\|\cdot\|)$,
	such that $A$ is convex, satisfies the property $UC^{*}$.
\end{proposition}

\begin{definition}\label{Definition:22}(\cite{Sintunavarat-Kumam})
	Let $A$ and $B$ be nonempty subsets of a metric space $(X,\rho)$ and $T:A\times A\to B$.
	A point $(x,y)\in A\times A$ is called a coupled best proximity point of $T$ in $A\times A$ if
	$\rho(x,T(x,y))=\rho(y, T(y,x))={\rm dist}(A, B)$.
\end{definition}

\begin{definition}\label{Definition:23}(\cite{Sintunavarat-Kumam})
	Let $A$ and $B$ be nonempty subsets of a metric space $(X,\rho)$. An ordered pair of maps $(F,G)$ is called an ordered pair of cyclic maps (or for short cyclic maps) if $F:A\times A\to B$ and $G:B\times B\to A$.
\end{definition}

\begin{definition}\label{Definition:24}(\cite{Gupta-Rajput-Kaurav,Sintunavarat-Kumam})
	Let $A$ and $B$ be nonempty subsets of a metric space $(X,\rho)$ an ordered pair of cyclic maps $(F,G)$
	is called an ordered pair of cyclic contraction maps if there exists $\alpha,\beta\in [0,1)$ with $\alpha+\beta<1$ so that the inequality 
	\begin{equation}\label{equation:1}
		\rho(F(x,y),G(u,v))\leq \alpha\rho(x,u)+\beta\rho(y,v)+(1-(\alpha+\beta)){\rm dist}(A,B)
	\end{equation}
	holds for every $(x,y)\in A\times A$ and $(u,v)\in B\times B$.
\end{definition}

The case when $\alpha=\beta$ is considered in \cite{Sintunavarat-Kumam}.

\begin{theorem}\label{Theorem:25}(\cite{Sintunavarat-Kumam})
	Let $A$ and $B$ be nonempty closed subsets of a complete metric space $(X,\rho)$,
	such that the ordered pairs $(A, B)$ and $(B, A)$ satisfy the property $UC^{*}$. 
	Let $F:A\times A\to B$, $G:B\times B\to A$ and $(F, G)$ be a cyclic contraction. 
	Then there exits coupled best proximity point $(x, y)$ of $F$ in $A\times A$ and 
	a coupled best proximity point $(u, v)$ of $G$ in $B\times B$ such that
	$\rho(x,u) + \rho(y,v) = 2{\rm dist}(A, B)$.
\end{theorem}

It seems that the properties $UC$ and $UC^{*}$ overlap for a wide class of sets with the key lemmas (Lemma \ref{Lemma:6} and Lemma \ref{Lemma:7}) from \cite{Eldred-Veeramani}.

By the imposed conditions in Theorem \ref{Theorem:14} and Theorem \ref{Theorem:25} it seems, at first glimpse, that there may be a gap in the proof 
of Theorem \ref{Theorem:14}, as far the generalization of Lemma \ref{Lemma:7}, in terms of a $UC^{*}$ property, is missing.
We will see in the  next section, that it is not the case, but the property $UC^{*}$ overlaps with the property $UC$,
when the sequences $\left\{x_n\right\}_{n=1}^\infty, \left\{z_n\right\}_{n=1}^\infty\subset A$ and $\left\{y_n\right\}_{n=1}^\infty\subset B$ are the iterated sequence, generated by 
a cyclic contraction $T$. We will show that we can remove the assumption $(A,B)$ to satisfy the property $UC^{*}$ in the conditions of Theorem \ref{Theorem:25}.

Actually the properties $UC$ and $UC^{*}$ can hold just for the two sets $A$ and $B$ that are the domain of the cyclic maps $T:A\cup B\to A\cup B$. We will try to find conditions which will ensure that whenever any ordered pair $(A,B)$ of subsets of a Banach space $(X,\|\cdot\|)$ satisfies the property $UC$, then $(X,\|\cdot\|)$ will be uniformly convex Banach space.

We will try to introduce a generalization of the notions of convexity, which will be different from $UC$, strict convexity or uniform convexity, but will insure existence and uniqueness of best proximity points for classes of cyclic maps $T:A\cup B\to A\cup B$.

\section{Connection between $UC$ and $UC^{*}$ properties}

We will start with some comments on the introduced in \cite{Sintunavarat-Kumam,Suzuki-Kikkawa-Vetro} properties $UC$ and $UC^{*}$.
The notions introduced in \cite{Sintunavarat-Kumam,Suzuki-Kikkawa-Vetro} search for some good properties to be satisfied by an ordered pair of sets $(A,B)$.  
The $UC$ states that for three sequences $\left\{x_n\right\}_{n=1}^\infty,\left\{z_n\right\}_{n=1}^\infty\subset A$, $\left\{y\right\}_{n=1}^\infty\subset B$, such that 
$\displaystyle\lim_{n\to\infty}\rho(x_n,y_n)=\lim_{n\to\infty}\rho(z_n,y_n)={\rm dist}(A,B)$, then
$\displaystyle\lim_{n\to\infty}\rho(x_n,z_n)=0$. 

The $UC^{*}$ requires, at first glimpse, a stronger property by
insisting the sequences $\left\{x_n\right\}_{n=1}^\infty,\left\{z_n\right\}_{n=1}^\infty\subset A$ to verify that for any $\varepsilon>0$ 
there is $N\in\mathbb{N}$ so that the inequality $\rho(x_m,z_n)\leq \varepsilon$
holds for all $m>n\geq N$. Actually the property $UC^{*}$ is included in the property $UC$ as we will see in  
Theorem \ref{Theorem:30}, because of the additional requirement on the sequences $\left\{x_n\right\}_{n=1}^\infty\subset A$, $\left\{y_n\right\}_{n=1}^\infty\subset B$ to satisfy
for every $\varepsilon>0$ there exists $N\in\mathbb{N}$, so that the inequality  
$\rho(x_m,y_n)\leq {\rm dist}(A,B)+\varepsilon$ holds for all $m>n\geq N$.

An auxiliary result in (\cite{Eldred-Veeramani}, Proposition 3.3) is that for a cyclic contraction map $T$, the iterated sequences
$\left\{T^{2n}x_0\right\}_{n=1}^\infty$, $\left\{T^{2n-1}x_0\right\}_{n=1}^\infty${,} for any arbitrary chosen initial guess point $x_0\in A\cup B${,}
are bounded ones. The authors in \cite{Eldred-Veeramani} apply Lemma \ref{Lemma:7} only for bounded sequences $\left\{x_n\right\}_{n=1}^\infty,\{z_n\}_{n=1}^\infty\subset A$, $\left\{y_n\right\}_{n=1}^\infty\subset B$ to show that the sequences $\left\{T^{2n}x_0\right\}_{n=1}^\infty$, $\left\{T^{2n-1}x_0\right\}_{n=1}^\infty$ are Cauchy ones.

A crucial lemma in \cite{Suzuki-Kikkawa-Vetro} presents a condition for a sequence to be a Cauchy one.

\begin{lemma}\label{Lemma:26}(\cite{Suzuki-Kikkawa-Vetro})
	Let $A$ and $B$ be subsets of a metric space $(X,\rho)$. Assume that the ordered pair $(A,B)$ has the property $UC$ and $\left\{x_n\right\}_{n=1}^\infty\subset A$, $\left\{y_n\right\}_{n=1}^\infty\subset B$, so that either one of following holds
	$$
	\lim_{m\to\infty}\sup_{n\geq m}\rho(x_m,y_n)={\rm dist}(A,B)\ \mbox{or}
	\ \lim_{n\to\infty}\sup_{m\geq n}\rho(x_m,y_n)={\rm dist}(A,B).
	$$
	Then $\{x_n\}$ is a Cauchy sequence. 
\end{lemma}

In the proof of the main result in \cite{Suzuki-Kikkawa-Vetro}, the authors show that
$$
\lim_{m\to\infty}\sup_{n\geq m}\rho\left(T^{2m}x,T^{2n+1}x\right)={\rm dist}(A,B),$$ 
which replaces the assumption the ordered
pair of sets $(A,B)$ to satisfies the $UC^{*}$ property.

As far as the application of Lemmas \ref{Lemma:6} and \ref{Lemma:7} or properties $UC$ or $UC^{*}$ are for 
the iterated sequences $\left\{T^{2n}x_0\right\}_{n=1}^\infty$, $\left\{T^{2n-1}x_0\right\}_{n=1}^\infty$, that 
are bounded ones, we will introduce a new property for an ordered pair of sets $(A,B)$, which will involve only bounded sequences $\left\{x_n\right\}_{n=1}^\infty, \left\{z_n\right\}_{n=1}^\infty\subset A$.

\begin{definition}\label{Definition:27}
	Let $A$ and $B$ be nonempty subsets of a metric space $(X,\rho)$. We say that the ordered pair $(A,B)$ 
	satisfies the bounded property $UC$ ($BUC$) if there holds one of the following
		\begin{itemize} 
			\item for every bounded sequences $\left\{x_n\right\}_{n=1}^\infty, \left\{z_n\right\}_{n=1}^\infty\subset A$ and
			an arbitrary sequence $\left\{y_n\right\}_{n=1}^\infty\subset B$, such that, whenever there holds 
			\begin{equation}\label{equation:2}
				\lim_{n\to\infty}\rho(x_n,y_n)=\lim_{n\to\infty}\rho(z_n,y_n)={\rm dist}(A,B),
			\end{equation} 
			then there holds $\lim_{n\to\infty}\rho(x_n,z_n)=0$
			\item there are no any bounded sequences $\left\{x_n\right\}_{n=1}^\infty, \left\{z_n\right\}_{n=1}^\infty\subset A$ to satisfy (\ref{equation:2})
		\end{itemize}.
\end{definition}

Let an ordered pair $(A,B)$ has the property $UC$. If equality (\ref{equation:2}) holds for bounded sequences $\left\{x_n\right\}_{n=1}^\infty, \left\{z_n\right\}_{n=1}^\infty\subset A$, then $(A,B)$  has the property $BUC$. 

If an ordered pair $(A,B)$ satisfies $BUC$, then it may happen that there are unbounded sequences $\left\{x_n\right\}_{n=1}^\infty, \left\{z_n\right\}_{n=1}^\infty\subset A$ so that the assumptions about property $UC$ are not satisfied and thus $(A,B)$ will have the $BUC$ property, but will not have the $UC$ property.

\begin{example}\label{Example:28}
	Let us consider the sets 
	$$
	A=\left\{(x,y)\in\mathbb{R}^2:y\geq \frac{1}{x-1},x>1\right\}\ \mbox{and}\ 
	B=\left\{(x,y)\in\mathbb{R}^2:y\geq \frac{1}{|x|},x<0\right\}.
	$$
	\begin{enumerate}
		\item If we consider $A,B\subset (\mathbb{R}^2,\|\cdot\|_2)$.
		Then the ordered pair $(A,B)$ has the $UC$ property.
		\item If we consider $A,B\subset (\mathbb{R}^2,\|\cdot\|_\infty)$.
		Then the ordered pair $(A,B)$ has not the $UC$ property.
	\end{enumerate}
\end{example}

In both cases of Example \ref{Example:28} the sequences that satisfy (\ref{equation:2}) are unbounded ones, therefore the ordered pair $(A,B)$ has the $BUC$ property. We can consider the sequences  $\{x_n\}_{n=1}^\infty=\{(\frac{1}{n}+1,n)\}_{n=1}^\infty\subset A$, $\{z_n\}_{n=1}^\infty=\{(\frac{1}{n}+1,n+1)\}_{n=1}^\infty\subset A$ and $\{y_n\}_{n=1}^\infty=\{(-\frac{1}{n},n)\}_{n=1}^\infty\subset A$.

First we will show that the conditions imposed on the sequences in the definition of the $UC^{*}$ property ensure 
that they are bounded sequences and therefore the set of the sequences in Definition \ref{Definition:27} is not the empty set 
for any ordered set $(A,B)$, satisfying the $UC^{*}$ property. 

\begin{lemma}\label{Lemma:29}
	Let $(X,\rho)$ be a metric space. If the sequences 
	$$
	\left\{x_n\right\}_{n=1}^\infty,\ \left\{z_n\right\}_{n=1}^\infty,\ \left\{y_n\right\}_{n=1}^\infty\subset X
	$$ 
	satisfy  
	\begin{equation}\label{equation:3}
		\lim_{n\to\infty}\rho(x_n,y_n)=a<\infty 
	\end{equation}
	\begin{equation}\label{equation:4}
		\lim_{n\to\infty}\sup_{n\leq k<m}\rho(z_m,y_k)=b<\infty, 
	\end{equation}
	then $\left\{x_n\right\}_{n=1}^\infty,\ \left\{z_n\right\}_{n=1}^\infty,\ \left\{y_n\right\}_{n=1}^\infty$ are bounded sequences.
\end{lemma}
\begin{proof}
By (\ref{equation:4}) it follows that there exists $n_0\in \mathbb{N}$ so that $\displaystyle\sup_{n_0<m}\rho(z_m,y_{n_0})<\infty$ and consequently the sequence $\{z_n\}_{n=1}^\infty$ is a bounded one. By similar arguments we get that $\{y_n\}_{n=1}^\infty$ is a bounded sequence too, From the boundedness of $\{y_n\}_{n=1}^\infty$ and (\ref{equation:3}) it follows that $\{x_n\}_{n=1}^\infty$ is a bounded sequence too.
\end{proof}

We have seen that, whenever an ordered pair $(A,B)$ satisfies the $UC$, then it satisfies $UC^{*}$.
	Therefore if an ordered pair $(A,B)$ satisfies the $BUC$ and $UC$ it will  satisfies $UC^{*}$ too.
	
	It is interesting what is the connection between $BUC$ and the conditions \ref{Definition:19}.\ref{enum-1}) and \ref{Definition:19}.\ref{enum-2}), without assuming the $UC$ property.

\begin{theorem}\label{Theorem:30}
	Let $(X,\rho)$ be a metric space and $A,B\subset X$. If the ordered pair $(A,B)$ satisfies the $BUC$ property, then the ordered pair $(A,B)$ satisfies the $RUC^*$ property
\end{theorem}

\begin{proof}
Let assume the contrary, i.e. $(A,B)$ satisfies the $BUC$, but does not satisfies the $UC^*$ property. 
Then there exist sequences $\{x_n\}_{n=1}^\infty,\{z_n\}_{n=1}^\infty\in A$ and $\{y_n\}_{n=1}^\infty\in B$, satisfying 
\begin{equation}\label{equation:5}
	\lim_{n\to\infty}\rho(x_n,y_n)={\rm dist}(A,B)\ \mbox{and} \ \lim_{n\to\infty}\sup_{n\leq k<m}\rho(z_m,y_k)={\rm dist}(A,B)
\end{equation}
and $\displaystyle\lim_{n\to\infty}\sup_{n\leq k<m}\rho(z_m,x_k)\neq0$. Therefore there exist subsequences 
$\{x_{k_i}\}_{i=1}^\infty\subset \{x_n\}_{n=1}^\infty$ and $\{z_{m_i}\}_{i=1}^\infty\subset \{z_n\}_{n=1}^\infty$,
so that 
\begin{equation}\label{equation:6}
	\lim_{i\to\infty}\rho(x_{n_i},z_{m_i})\neq 0.
\end{equation}
Let us consider the sequences $\{x_{k_i}\}_{i=1}^\infty, \{z_{m_j}\}_{j=1}^\infty\subset A$ 
and $\{y_{k_i}\}_{i=1}^\infty\subset B$.
By (\ref{equation:5}) it follows that
\begin{equation}\label{equation:7}
	\lim_{n\to\infty}\rho(x_{k_i},y_{k_i})={\rm dist}(A,B)
\end{equation}
and
for every $\varepsilon>0$, there is $N\in\mathbb{N}$, so that for all $i,j\geq N$ there holds the inequality
\begin{equation}\label{equation:8}
	\rho(z_{m_j},y_{k_i})\leq{\rm dist}(A,B)+\varepsilon.
\end{equation}
From (\ref{equation:8}) we get that $\lim_{i\to\infty}\rho(z_{m_i},y_{k_i})={\rm dist}(A,B)$.
By Lemma \ref{Lemma:29}, (\ref{equation:7}) and (\ref{equation:8}) it follows that the sequences are bounded ones.
From the assumption of the theorem that the ordered pair $(A,B)$ has the $BUC$ property and the sequences
$\{x_{k_i}\}_{i=1}^\infty, \{z_{m_j}\}_{j=1}^\infty\subset A$ and $\{y_{k_i}\}_{i=1}^\infty\subset B$ are bounded ones
it follows that 
$\lim_{i\to\infty}\rho(x_{n_i},z_{m_i})=0$, which is a contradiction with (\ref{equation:7}).
\end{proof}

Having in mind Lemma \ref{Lemma:29}, the conditions imposed on the sequences in the definition of the $UC^{*}$ property ensure 
that they are bounded ones and thus it follows that from the $UC$ property {it} follows the $UC^{*}$ property, that is $UC^{*}$ is included in the $UC$ property. Consequently we can replace in Theorem \ref{Theorem:25} the assumption that $(A,B)$ satisfies the $UC^{*}$ by $(A,B)$ satisfies the $UC$.

From Theorem \ref{Theorem:30} and \cite{Eldred-Veeramani}, where the authors have proven that the iterated sequence $\{T^nx\}_{n=1}^\infty$, for every an arbitrary chosen initial guess $x\in A\cup B$, is a bounded one for the maps investigated in Theorem \ref{Theorem:5}, it follows that we can present a generalization of their result. We will illustrate in the last section with an example that the next theorem is actually a generalization of Theorem \ref{Theorem:5}. We will present an ordered pair $(A,B)$ and a cyclic map $T$, so that $(A,B)$ has the $BUC$ property, but has not the $UC$ property and $T$ satisfies the contracive condition in Theorem \ref{Theorem:5}, i.e. (\ref{equation:9}). 

\begin{theorem}\label{Theorem:31}
	Let $A$ and $B$ be nonempty closed subsets of a complete metric space $(X,\rho)$,
	such that the ordered pair $(A, B)$ satisfies the property $BUC$. 
	Let $T:A\cup B\to A\cup B$ be a cyclic map and there exists $k\in [0,1)$, so that the inequality
	\begin{equation}\label{equation:9}
		\rho(Tx,Ty)\leq k\rho(x,y)+(1-k){\rm dist}(A,B)
	\end{equation}
	holds for all $x\in A$ and $y\in B$. Then there is a unique best proximity point $x$ of $T$ in $A$, the sequence of successive iterations $\{T^{2n}x_0\}_{n=1}^\infty$ converges to $x$ for any initial guess $x_0\in A$. There is at least one best proximity point in $B$ of $T$ and if the ordered pair $(B,A)$ has the $BUC$ property, then this point is unique.
\end{theorem}

We will show that the iterated sequence $\{T^nx\}_{n=1}^\infty$, for every arbitrary chosen initial guess $x\in A\cup B$ is a bounded one for the maps investigated by \cite{Suzuki-Kikkawa-Vetro} in Theorem \ref{Theorem:14}.

\begin{lemma}\label{Lemma:32}
	Let $(X,\rho )$ be a metric space and let $A$ and $B$ be nonempty subsets of $X$ and $T:A\cup B\to A\cup B$ be a cyclic map, such that there is $\lambda\in [0,1)$ such that the inequality
	\begin{equation}\label{equation:10}
		\rho(Tx,Ty)\leq\lambda\max\{\rho(x,y),\rho(x,Tx),\rho(y,Ty)\}+(1-\lambda){\rm dist}(A,B)
	\end{equation}
	holds for every $x\in A$ and $y\in B$.
	Then the iterated sequence $\{T^nx\}_{n=1}^\infty$ is a bounded one for every $x\in A \cup B$.
\end{lemma}

\begin{proof}
Just to fit some of the formulas into the text field we will use the notation $d={\rm dist}(A,B)$.
Let $x\in A \cup B$ be an arbitrary chosen point. Let us denote $x_n=T^nx$, for $n=0,1,2,\dots$. 

We will show that 
\begin{equation}\label{equation:11}
	\sup_{n\in \mathbb{N}}\rho(x_n,x_{n-1})=\rho(x_{1},x_{0})<\infty.
\end{equation}

For every $n\in \mathbb{N}$ by (\ref{equation:10}) we have the inequality 
\begin{equation}\label{equation:12}
	\rho(x_n,x_{n-1})\leq\lambda\max\{\rho(x_{n-1},x_{n-2}),\rho(x_{n-1},x_n),\rho(x_{n-1},x_{n-2})\}+(1-\lambda)d,
\end{equation}
where $d$ stands for ${\rm dist}(A,B)$.
The inequality (\ref{equation:12}) reduces to one of the inequalities, either
\begin{equation}\label{equation:13}
	\rho(x_n,x_{n-1})\leq\lambda\rho(x_{n-1},x_n)+(1-\lambda){\rm dist}(A,B)
\end{equation}
or 
\begin{equation}\label{equation:14}
	\rho(x_n,x_{n-1})\leq\lambda\rho(x_{n-1},x_{n-2})+(1-\lambda){\rm dist}(A,B).
\end{equation}
If there holds (\ref{equation:13}), then 
\begin{equation}\label{equation:15}
	\rho(x_n,x_{n-1})={\rm dist}(A,B).
\end{equation}

For every arbitrary chosen $n\in\mathbb{N}$ there are two cases: for some $m<n$ there holds
either
\begin{equation}\label{equation:16}
	\rho(x_m,x_{m-1})={\rm dist}(A,B)
\end{equation}
or
\begin{equation}\label{equation:17}
	\rho(x_m,x_{m-1})\leq\lambda\rho(x_{m-1},x_{m-2})+(1-\lambda){\rm dist}(A,B).
\end{equation}
If there holds (\ref{equation:17}) for every $m<n$, then we get the chain of inequalities
$$
\begin{array}{rl}
	\rho(x_n,x_{n-1})\leq&\lambda\rho(x_{n-1},x_{n-2})+(1-\lambda){\rm dist}(A,B)\\
	\leq&\lambda^2\rho(x_{n-2},x_{n-3})+(1-\lambda^2){\rm dist}(A,B)\\
	\cdots&\cdots\\
	\leq&\lambda^{k}\rho(x_{n-k},x_{n-k+1})+(1-\lambda^{k}){\rm dist}(A,B)\\
	\cdots&\cdots\\
	\leq&\lambda^{n-2}\rho(x_{2},x_{1})+(1-\lambda^{n-2}){\rm dist}(A,B)\\
	\leq&\lambda^{n-1}\rho(x_{1},x_{0})+(1-\lambda^{n-1}){\rm dist}(A,B)\leq\rho(x_{1},x_{0}).
\end{array}
$$
If for some $m< n$ the case (\ref{equation:16}) is valid, then we get
$$
\begin{array}{rl}
	\rho(x_n,x_{n-1})\leq&\lambda\rho(x_{n-1},x_{n-2})+(1-\lambda){\rm dist}(A,B)\\
	\leq&\lambda^2\rho(x_{n-2},x_{n-3})+(1-\lambda^2){\rm dist}(A,B)\\
	\cdots&\cdots\\
	\leq&\lambda^{m-n}\rho(x_{m},x_{m-1})+(1-\lambda^{m-n}){\rm dist}(A,B)\\
	\leq&{\rm dist}(A,B)\leq\rho(x_{1},x_{0}),
\end{array}
$$
which finishes the proof of (\ref{equation:11}).

For every $n$ from (\ref{equation:10}) we have
$$
\rho(x_{2n},x_1)\leq\lambda\max\{\rho(x_{2n-1},x_0),\rho(x_{2n-1},x_{2n}),\rho(x_1,x_0)\}+(1-\lambda){\rm dist}(A,B).
$$
By (\ref{equation:11}) we get
$\rho(x_{2n},x_1)\leq\lambda\max\{\rho(x_{2n-1},x_0),\rho(x_1,x_0)\}+(1-\lambda){\rm dist}(A,B)$.
From the triangle inequality and (\ref{equation:11}) it follows
$$
\rho(x_{2n-1},x_0)\leq \rho(x_{2n-1},x_{2n})+\rho(x_{2n},x_1)+\rho(x_{1},x_0)\leq 2\rho(x_{1},x_0)+\rho(x_{2n},x_1).
$$
Therefore we get $\rho(x_{2n},x_1)\leq (2\lambda+1)\rho(x_1,x_0)+\lambda\rho(x_{2n},x_1)+(1-\lambda){\rm dist}(A,B)$, 
i.e.
$$
\rho(x_{2n},x_1)\leq \frac{2\lambda+1}{1-\lambda}\rho(x_1,x_0)+{\rm dist}(A,B).
$$ 
Consequently $\left\{x_{2n}\right\}_{n=0}^\infty$ is a bounded sequence. 

By similar arguments we can prove that $\{x_{2n+1}\}_{n=0}^\infty$ is a bounded sequence too.
\end{proof}

From Theorem \ref{Theorem:30} and Lemma \ref{Lemma:32} we can present a 
generalization of Theorem \ref{Theorem:14}.

\begin{theorem}\label{Theorem:33}(\cite{Suzuki-Kikkawa-Vetro}) 
	Let $A$ and $B$ be nonempty closed subsets of a complete metric space $(X,\rho)$,
	such that the ordered pair $(A, B)$ satisfies the property $BUC$. 
	Let $T:A\cup B\to A\cup B$ be a cyclic map and there exists $k\in [0,1)$, so that the inequality
	$$
	\rho(Tx,Ty)\leq k\max\{\rho(x,y),\rho(x,Tx),\rho(y,Ty)\}+(1-k){\rm dist}(A,B)
	$$
	holds for all $x\in A$ and $y\in B$. Then there is a unique best proximity point $x$ of $T$ in $A$, the sequence of successive iterations $\{T^{2n}x_0\}_{n=1}^\infty$ converges to $x$ for any initial guess $x_0\in A$. There is at least one best proximity point in $B$ of $T$ and if the ordered pair $(B,A)$ has the $BUC$ property, then this point is unique.
\end{theorem}

We will show that the iterated sequences generated by the maps in Theorem \ref{Theorem:25} are bounded ones.
Following a smart idea, from \cite{Petrusel}, that connects coupled fixed points and fixed points, we will try to apply
the technique from \cite{Petrusel} to show that the maps involved in Theorem \ref{Theorem:25} can be considered as like as the maps 
investigated in \cite{Eldred-Veeramani}.

Let us point out that in instead of considering two maps $F:A\times A\to B$ and $G:B\times B\to A$ (Definition \ref{Definition:24}, Theorem \ref{Theorem:25} we can consider just one map $f:(A\times A)\cup (B\times B)\to (A\times A)\cup (B\times B)$ defined by
	$$
	f(x,y)=\left\{
	\begin{array}{rl}
		F(x,y),&x,y\in A\\
		G(x,y),&x,y\in B.
	\end{array}
	\right.
	$$

\begin{lemma}\label{Lemma:34}
	Let $(X,\rho )$ be a metric space and let $A$ and $B$ be nonempty subsets of $X$. 
	Let $f$ be a cyclic contraction, satisfying (\ref{equation:1}), i.e.
	\begin{equation}\label{equation:18}
		\rho(f(x,x'),f(y,y'))\leq\alpha\rho(x,y)+\beta\rho(x',y')+(1-(\alpha+\beta)){\rm dist}(A,B)
	\end{equation}
	for every $x,x'\in A$, $y,y'\in B$ and some $\alpha,\beta\geq 0$ and $\alpha+\beta\in[0,1)$.
	Then
	the iterated sequences $\{x_n\}_{n=1}^\infty$ and $\{y_n\}_{n=1}^\infty$ are bounded ones.
\end{lemma}

\begin{proof} Let us consider the metric space $(X\times X, d)$, where $d((x,y),(u,v))=\rho(x,u)+\rho(y,v)$.
Let us define the map $T:X\times X\to X\times X$, by $T(x,y)=(f(x,y),f(y,x))$. The map $T$ is a cyclic map as
$T(A\times A)\subseteq B\times B$ and $T(B\times B)\subseteq A\times A$.
There holds
\begin{equation}\label{equation:19}
	\begin{array}{lll}
		{\rm dist}(A\times A,B\times B)&=&\inf\{d((a,a^\prime),(b,b^\prime)):a,a^\prime\in A,b,b^\prime\in B\}\\
		&=&\inf\{\rho(a,b)+\rho(a^\prime,b^\prime):a,a^\prime\in A,b,b^\prime\in B\}\\
		&=&2\inf\{\rho(a,b):a,\in A,b\in B\}
		=2{\rm dist}(A,B).
	\end{array}
\end{equation}
For any two arbitrary chosen $x=(x,x^\prime)\in A\times A$ and $y=(y,y^\prime)\in B\times B$, by using 
(\ref{equation:18}) and (\ref{equation:19}), there holds the chain of inequalities
$$
\begin{array}{lll}
	d(Tx,Ty)&=&d((f(x,x^\prime),(f(x^\prime,x)),(f(y,y^\prime),(f(y^\prime,y)))\\
	&=&\rho(f(x,x^\prime),f(y,y^\prime))+\rho(f(x^\prime,x),f(y^\prime,y))\\
	&\leq&\alpha\rho(x,y)+\beta\rho(x',y')+(1-(\alpha+\beta)){\rm dist}(A,B)\\
	&&+\alpha\rho(x',y')+\beta\rho(x,y)+(1-(\alpha+\beta)){\rm dist}(A,B)\\
	&=&(\alpha+\beta)(\rho(x,y)+\rho(x',y'))+2(1-(\alpha+\beta)){\rm dist}(A,B)\\
	&=&(\alpha+\beta)d((x,x^\prime),(y,y^\prime))+(1-(\alpha+\beta)){\rm dist}(A\times A,B\times B).
\end{array}
$$
Consequently the cyclic map $T$ satisfies the conditions imposed in \cite{Eldred-Veeramani} and
according to \cite{Eldred-Veeramani} the iterated sequence 
$$
u_n=(x_n,y_n)=Tu_{n-1}=T(x_{n-1},y_{n-1})=(f(x_{n-1},y_{n-1}),f(y_{n-1},x_{n-1}))
$$ 
for every arbitrary chosen $u_0=(x_0,y_0)\in A\times A$ is a bounded sequence.
\end{proof}

From Theorem \ref{Theorem:30} and Lemma \ref{Lemma:34} we can present a 
generalization of Theorem \ref{Theorem:25}.

\begin{theorem}\label{Theorem:35}(\cite{Suzuki-Kikkawa-Vetro}) 
	Let $A$ and $B$ be nonempty closed subsets of a complete metric space $(X,\rho)$,
	such that the ordered pairs $(A, B)$ and $(B, A)$ satisfy the property $BUC$. 
	Let $F:A\times A\to B$, $G:B\times B\to A$ and $(F, G)$ be a cyclic contraction. 
	Then there exits a coupled best proximity point $(x, y)$ of $F$ in $A\times A$ and 
	a coupled best proximity point $(u, v)$ of $G$ in $B\times B${\color{blue},} such that $\rho(x,u) + \rho(y,v) = 2{\rm dist}(A, B)$
\end{theorem}

\section{Characterization of the $UC$ property}

We have seen in the examples, that it is possible to have the underlying Banach space $(X,\|\cdot\|)$ not to be uniformly convex, but some particularly chosen ordered pair of subsets $(A,B)$ to be either $UC$ one or $BUC$ one. We will try to find some sufficient conditions that will ensure an ordered pair of subsets $(A,B)$ to be a $UC$ one in an arbitrary Banach space. 

\begin{definition}\label{Definition:36}
	Let $(X,\| \cdot \|)$ be a Banach space and $A\subset X$. We say that $A$ is uniformly convex set if for every $\epsilon>0$ there exists $\eta(\epsilon)>0$, such that for every $x,y \in A$, satisfying the inequality $\|x-y\|\geq \epsilon$ there holds $B\left(\frac{x+y}{2},\eta(\epsilon)\right)\subset A$.
\end{definition}

From Definition \ref{Definition:1} and Definition \ref{Definition:36} it is easy to observe that for every $\varepsilon>0$ and for any $x,y\in B_{(X,\|\cdot\|)}$, such that $\|x-y\|\geq \varepsilon$, there exists $\eta(\varepsilon)>0$, so that $B\left(\frac{x+y}{2},\eta(\epsilon)\right)\subset B_X$.
We can choose $\eta(\varepsilon)=\delta_X(\varepsilon)$. Thus the unit ball in a uniformly convex Banach space satisfies Definition \ref{Definition:36}.

We have seen in Example \ref{Example:16}, that if the sets $A$ or $B$ have some good geometric properties then the ordered pair $(A,B)$ will be an $UC$ one.

\begin{theorem}\label{Theorem:37}
	Let $(X,\|\cdot\|)$ be a Banach space, $A,B\subset X$ and $A$ be a uniformly convex set. Then the ordered pair $(A,B)$ has the $UC$ property.
\end{theorem}

\begin{proof}
Let us assume the contrary, i.e. there are sequences $\{x_n\}_{n=1}^\infty\subset A$, $\{z_n\}_{n=1}^\infty\subset A$ and $\{y_n\}_{n=1}^\infty\subset B$, such that 
$$
\lim_{n\to\infty}\|x_n-y_n\|=\lim_{n\to\infty}\|z_n-y_n\|={\rm dist}(A,B),
$$ 
but the sequence $\{\|x_n-z_n\|\}_{n=1}^\infty$ does not converge to zero.
Then there exists $\varepsilon_0>0$, so that for every $N\in\mathbb{N}$ there is $n>N$ and the inequality 
$\left\|x_n-z_n\right\|>\varepsilon_0$ holds true. From the assumption, that $A$ is a uniformly convex set it follows that for $\varepsilon_0$, there is $\eta(\varepsilon_0)$ so that the inclusion $B\left(\frac{x_n+z_n}{2},\eta(\varepsilon_0)\right)\subset A$ holds. Consequently we can write the chain of inequalities 
\begin{equation}\label{equation:20}
	\begin{array}{lll}
		\eta(\epsilon_0)+{\rm dist}(A,B)&\leq& {\rm dist}\left(\frac{x_n+z_n}{2},B\right)\leq\left\|\displaystyle\frac{x_n+z_n}{2}-y_n\right\|\\[10pt]
		&\leq& \displaystyle\frac{\|x_n-y_n\|+\|z_n-y_n\|}{2}.
	\end{array}
\end{equation}
By the choice of the sequences $\{x_n\}_{n=1}^\infty$, $\{z_n\}_{n=1}^\infty$ and $\{y_n\}_{n=1}^\infty$ it follows that there exists $N_1(\varepsilon_0)\in \mathbb{N}$, such that for every $n>N_1$ there holds the inequality
\begin{equation}\label{equation:21}
	\frac{\|x_n-y_n\|+\|z_n-y_n\|}{2}<\eta(\varepsilon_0)+{\rm dist}(A,B),
\end{equation}
which is a contradiction with (\ref{equation:20}). Thus $\lim_{n\to\infty}\|x_n-z_n\|=0$ and we get that the ordered pair $(A,B)$ has the $UC$ property.
\end{proof}

By the results up to now it follows that if $A$ is a uniformly convex set, then every ordered pair of sets $(A,B)$ has the $BUC$ property and the $UC$ property. It seems that the assumption $A$ to be a uniformly convex set is too restrictive. We will try to find a weaker property, which will ensure that the ordered pair of sets $(A,B)$ has the $BUC$ property, without being $UC$ ordered pair.

\begin{definition}\label{Definition:38}
	Let $(X,\| \cdot \|)$ be a Banach space and $A\subset X$. We say that a function $\phi$ has the positive property about the set $A$ if $\phi: A\times \mathbb{R}^+ \to \mathbb{R}^+$ is such that for every bounded subset $A'\subset A$ and every $\epsilon_0>0$ there holds the inequality 
	$$
	\displaystyle\inf\{\phi(x,\epsilon):x\in A',\epsilon\geq \epsilon_0\}>0.
	$$
\end{definition}

We will illustrate Definition \ref{Definition:38} by an example.

\begin{example}\label{Example:39}
	Let us consider the space $X=(\mathbb{R}^2,\|\cdot\|_\infty)$,
	the set 
	$$
	\displaystyle A=\left\{(x,y)\in\mathbb{R}^2:\frac{1}{x}\leq y,x>0\right\}
	$$
	and the function $\displaystyle\phi:A\times\mathbb{R}^+  \rightarrow \mathbb{R}^+$ be defined by 
	$\displaystyle\phi(x,\varepsilon)=\frac{\varepsilon^2}{320+5\varepsilon^2+5\|x\|^3}$.
	Then $\phi$ has the positive property about $A$.
\end{example}

Let us define the function $\phi_1:\mathbb{R}^+\times\mathbb{R}^+\rightarrow \mathbb{R}^+$ by 
$\displaystyle\phi_1(r,\varepsilon)=\frac{\varepsilon^2}{320+5\varepsilon^2+5r^3}$. Let us consider the function
\begin{equation}\label{equation:22}
	\phi(x,\varepsilon)=\phi_1(\|x\|,\varepsilon).
\end{equation}
From the boundedness of $A'\in A$ it follows the existence of $r_1=r_1(A^\prime)$, so that 
$$
\displaystyle\sup\{\|x\|:x\in A'\}=r_1<\infty.
$$ 
For every $\varepsilon_1>0$ there holds 
$$\displaystyle\inf\{\phi_1(r,\varepsilon):0\leq r\leq r_1,0<\varepsilon_1\leq\varepsilon\}=\phi_1(r_1,\varepsilon_1)>0.$$
By (\ref{equation:22}) and the definition of $r_1$ it follows that 
$$\displaystyle\inf\{\phi(x,\varepsilon):x\in A',0<\varepsilon_1\leq\varepsilon\}\geq\phi_1(\sup\{\|x\|:x\in A'\},\varepsilon_1)>0.$$
Consequently $\phi$ has positive property about $A$.

\begin{definition}\label{Definition:40}
	Let $(X,\| \cdot \|)$ be a Banach space, $A\subset X$ and $\phi: A\times \mathbb{R}^+ \to \mathbb{R}^+$ has the positive property about $A$. We say that $A$ is a uniformly convex set about $\phi$, if for every $\varepsilon>0$ and every $x,y \in A$, satisfying $\|x-y\|\geq \varepsilon$, there holds $B\left(\frac{x+y}{2},\phi(\frac{x+y}{2},\varepsilon)\right)\subset A$.
\end{definition}

We will show that the set in Example \ref{Example:39} is uniformly convex about $\phi$.

\begin{example}\label{Example:41}
	Let us consider the space $X=(\mathbb{R}^2,\|\cdot\|_\infty)$, the set 
	$$
	\displaystyle A=\left\{(x,y)\in\mathbb{R}^2: \frac{1}{x}\leq y, x>0 \right\}
	$$
	and the function $\phi:A\times\mathbb{R}^+  \to \mathbb{R}^+$ be defined by 
	$\displaystyle\phi(x,\varepsilon)=\frac{\varepsilon^2}{320+5\varepsilon^2+5\|x\|^3}$.
	Then $A$ is uniformly convex about $\phi$.
	
\end{example}

We have proven in Example \ref{Example:39} that $\phi$ has the positive property about $A$. It remains to show that Then $A$ is uniformly convex about $\phi$.

Let $\varepsilon>0$ be arbitrary chosen and let $p_1=(x_1,y_1)\in A$, $p_2=(x_2,y_2)\in A$ satisfy $\|p_1-p_2\|=\varepsilon$. Let us put $p_3=(x,y)=\frac{p_1+p_2}{2}=\left(\frac{x_1+x_2}{2},\frac{y_1+y_2}{2}\right)$. 

There are two cases: either 1) $|x_1-x_2|=\varepsilon$ or 2) $|y_1-y_2|=\varepsilon$ is fulfilled. 

Let us consider first the case 1) $|x_1-x_2|=\varepsilon$. 

Without loss generality we can assume that $x_1<x_2$. Then for $x$ we obtain the equalities 
\begin{equation}\label{equation:23}
	\begin{array}{lll}
		x_2&=&x+\frac{\varepsilon}{2}\\
		x_1&=&x-\frac{\varepsilon}{2}.
	\end{array}
\end{equation} 
From $p_1\in A$ it follows that $x_1>0$ and therefore $\displaystyle0<x-\frac{\varepsilon}{2}$. Thus $\varepsilon<2x$. 

From the inequalities $\displaystyle y_1\geq \frac{1}{x_1}$, $\displaystyle y_2\geq \frac{1}{x_2}$ and by using (\ref{equation:23}) we get that there hold
$$
\displaystyle y_1\geq \frac{1}{x-\varepsilon},\ \displaystyle y_2\geq \frac{1}{x+\varepsilon}\ \mbox{and}\ \displaystyle y\geq \frac{x}{x^2-\frac{\varepsilon^2}{4}}.
$$ 
Thus if $\displaystyle p_3=(x,y)$, then $y\geq \frac{x}{x^2-\frac{\varepsilon^2}{4}}$ and $0<\varepsilon<2x$.

If there holds the case 2) $|y_1-y_2|=\varepsilon$, by similar arguments, we get that there hold 
$x\geq \frac{y}{y^2-\frac{\varepsilon^2}{4}}$ and $0<\varepsilon<2y$, provided that $\displaystyle p_3=(x,y)$.

Consequently for any $\varepsilon>0$ and any $\displaystyle p_1\in A$, $p_2\in A$, such that $\|p_1-p_2\|=\varepsilon$ we get that $p_3=\frac{p_1+p_2}{2}\in P(\varepsilon)=C(\varepsilon)\cup D(\varepsilon)$, where 
$$
C(\varepsilon)=\left\{(x,y)\in\mathbb{R}^2:x\geq \frac{y}{y^2-\frac{\varepsilon^2}{4}},0<\varepsilon<2y\right\}
$$
and
$$
D(\varepsilon)=\left\{(x,y)\in\mathbb{R}^2:y\geq \frac{x}{x^2-\frac{\varepsilon^2}{4}},0<\varepsilon<2x\right\}.
$$

Thus to show that $A$ is uniformly convex about $\phi$, it is enough to estimate $\phi(p,\varepsilon)$, for $p\in P(\varepsilon)$.

For every $\varepsilon>0$ and $p\in P(\varepsilon)$ there holds the inequality  $\displaystyle\phi(p,\varepsilon)\leq {\rm dist}\left(p,\overline{A}\right)$, where $\overline{A}=\left\{(x,y)\in\mathbb{R}^2: \frac{1}{x}=y, x>0 \right\}$ and thus $B(p,\phi(p,\varepsilon))\in A$. 

\begin{theorem}\label{Theorem:42}
	Let $(X,\|\cdot\|)$ be a Banach space, $A,B\subset X$, $\phi: A\times \mathbb{R}^+ \rightarrow \mathbb{R}^+$ has the positive property  about $A$ and $A$ be a uniformly convex set about $\phi$. Then the ordered pair $(A,B)$ satisfies the $BUC$ property.
\end{theorem}

\begin{proof} We will prove the theorem by assuming the contrary, i.e. let the ordered pair $(A,B)$ does not satisfy the $BUC$ property and $A$ be uniformly convex set about $\phi$. Then there are three bounded sequences $\left\{x_n\right\}_{n=1}^\infty, \left\{z_n\right\}_{n=1}^\infty\subset A$ and $\left\{y_n\right\}_{n=1}^\infty\subset B$ such that
\begin{equation}\label{equation:24}
	\lim_{n\to\infty}\|x_n-y_n\|=\lim_{n\to\infty}\|z_n-y_n\|={\rm dist}(A,B),
\end{equation}
but the sequence $\{\|x_n-z_n\|\}_{n=1}^\infty$ does not converge to zero.

Therefore there exists $\varepsilon_0>0$ such that for every $N\in\mathbb{N}$ there is $n>N$ so that the inequality 
$\left\|x_n-z_n\right\|\geq\epsilon_0$ holds true. Thus there is a subsequence $\{\|x_{n_k}-z_{n_k}\|\}_{k=1}^\infty$ so that $\|x_{n_k}-z_{n_k}\|\geq\varepsilon_0$ for all $k\in\mathbb{N}$.  
The sequence $\left\{\frac{x_{n_k}+z_{n_k}}{2}\right\}_{k=1}^\infty\in A$ is a bounded sequence too. 
By the assumption that $\phi$ has the positive property about $A$ and that $\left\{\frac{x_{n_k}+z_{n_k}}{2}\right\}_{k=1}^\infty\in A$ is a bounded sequence it follows that there exists $\delta_0$, such that the inequality
\begin{equation}\label{equation:25}
	\inf\left\{\phi(a,\varepsilon):a\in\left\{\frac{x_n+z_n}{2}\right\}_{n=1}^\infty,\varepsilon\geq\varepsilon_0\right\}=\delta_0>0
\end{equation}
holds true.

From the assumption that $A$ is uniformly convex set about $\phi$ and the fact that $\left\{\frac{x_{n_k}+z_{n_k}}{2}\right\}_{k=1}^\infty\in A$ is a bounded sequence, satisfying the inequality $\|x_{n_k}-z_{n_k}\|\geq\varepsilon_0$ we get the inclusions
$$
B\left(\frac{x_n+z_n}{2},\delta_0\right)\subseteq B\left(\frac{x_n+z_n}{2},\phi\left(\frac{x_n+z_n}{2},\varepsilon_0\right)\right)\subset A.
$$

Consequently for all $k\geq N_1$ there holds the chain of inequalities 
	\begin{equation}\label{equation:26}
		\begin{array}{lll}
			\delta_0+{\rm dist}(A,B)&\leq& {\rm dist}\displaystyle\left(\frac{x_{n_k}+z_{n_k}}{2},B\right)\leq\left\|\frac{x_{n_k}+z_{n_k}}{2}-y_n\right\|\\[10pt]
			&\leq& \displaystyle\frac{\|x_{n_k}-y_{n_k}\|+\|z_{n_k}-y_{n_k}\|}{2}.
		\end{array}
	\end{equation}
	By (\ref{equation:24}) it follows that for $\delta_0$ there exists $N_1\in \mathbb{N}$ such that for every $k\geq N_1$ the inequality
	$$
	\frac{\|x_{n_k}-y_{n_k}\|+\|z_{n_k}-y_{n_k}\|}{2}<\delta_0+{\rm dist}(A,B)
	$$
	holds true.

Which is a contradiction with (\ref{equation:26}) and thus the sequence $\{\|x_n-z_n\|\}_{n=1}^\infty$ converges to zero.
\end{proof}

\begin{example}\label{Example:43}
	Let us consider the space $X=(\mathbb{R}^2,\|\cdot\|_\infty)$ and the sets 
	$\displaystyle A=\left\{(x,y)\in\mathbb{R}^2: \frac{1}{x}\leq y, x>0 \right\}$, $\displaystyle B=\left\{(x,y)\in\mathbb{R}^2: \frac{1}{x+1}-1\geq y, x>-1\right\}$,
	then the ordered pair $(A,B)$ has the $BUC$ property, but does not have the $UC$ property
\end{example}
We have proven in Examples \ref{Example:39} and \ref{Example:41} that there exists a function $\phi$, so that $\phi$ has the positive property about $A$
and $A$ is uniformly convex about $\phi$. By Theorem \ref{Theorem:42} it follows that the ordered pair $(A,B)$ has the $BUC$ property, i.e.
for every two bounded sequences $\{x_n\}_{n=1}^\infty,\{z_n\}_{n=1}^\infty\subset A$ and a sequence $\{y_n\}_{n=1}^\infty\subset B$, such that $\lim_{n\to\infty}\|x_n-y_n\|=\lim_{n\to\infty}\|z_n-y_n\|={\rm dist}(A,B)$
the sequence $\{\|x_n-z_n\|\}_{n=1}^\infty$ converges to zero.

It remains to show that the ordered pair $(A,B)$ does not have the $UC$ property, i.e. we will show that there are unbounded sequences 
$\{x_n\}_{n=1}^\infty\subset A$, $\{z_n\}_{n=1}^\infty\subset A$ and a sequence $\{y_n\}_{n=1}^\infty\subset B$, such that 
$$
\lim_{n\to\infty}\|x_n-y_n\|=\lim_{n\to\infty}\|z_n-y_n\|={\rm dist}(A,B),
$$ 
but the sequence $\{\|x_n-z_n\|\}_{n=1}^\infty$ does not converge to zero.

Let consider the sequences $\{a_n\}_{n=0}^\infty,\ \{b_n\}_{n=0}^\infty$ and $\{c_n\}_{n=0}^\infty$, defined as follows 
$$\textstyle a_n=\left(n\ ,\ \frac{1}{n}\right)\ ,\ b_n=\left(n+1\ ,\ \frac{1}{n+1}\right)\in A\ ,\ c_n=\left(n\ ,\ \frac{1}{n+1}-1\right)\in B,\ \mbox{for}\ n\in\mathbb{N}.$$

From $\|a_n-c_n\|_\infty=\left\|\left(0\ ,\ \frac{1}{n}-\frac{1}{n+1}+1\right)\right\|_\infty=\frac{1}{n}-\frac{1}{n+1}+1$ and $\|b_n-c_n\|=\left\|\left(1\ ,\ 1\right)\right\|_\infty=1$ it follows that $\displaystyle \lim_{n\to \infty}\|a_n-c_n\|_\infty=1$ and $\displaystyle \lim_{n\to \infty}\|b_n-c_n\|_\infty=1$. 

By $\|b_n-a_n\|_\infty=\left\|\left(1\ ,\ \frac{1}{n+1}-\frac{1}{n}\right)\right\|_\infty\geq 1$ we get $\lim_{n\to \infty}\|b_n-a_n\|_\infty\neq 0$. 

From ${\rm dist}(A,B)=1$ it follows that $\displaystyle\lim_{n\to \infty}\|a_n-c_n\|_\infty={\rm dist}(A,B)$, $\displaystyle\lim_{n\to \infty}\|b_n-c_n\|_\infty={\rm dist}(A,B)$ and $\displaystyle\lim_{n\to \infty}\|b_n-a_n\|_\infty\neq 0$ and consequently the ordered pair $(A,B)$ does not have the $UC$ property.

\section{The property $UC$ and uniform convexity of the underlying space}

We will show that in some cases the validity of the $UC$ property in a Banach space $(X,\|\cdot\|)$ leads to a conclusion that $(X,\|\cdot\|)$ is a uniformly convex Banach space.
We will start with two technical lemmas.

\begin{lemma}\label{lem:4}
	Let $(X,\|\cdot\|)$ be a Banach space. Let $\{x_n\}_{n=1}^\infty,\{y_n\}_{n=1}^\infty \in X$ be sequences, such that for every $n\in\mathbb{N}$ there hold $\|x_n\|\leq a$, $\|y_n\|\leq b$ and $\lim_{n\rightarrow \infty}\|x_n+y_n\|=a+b$, where $a,b\in [0,+\infty)$.
	Then $\lim_{n\rightarrow \infty}\|x_n\|=a$ and $\lim_{n\rightarrow \infty}\|y_n\|=b$.
\end{lemma}
\begin{proof}
From the chain of inequalities
$$a=\lim_{n\to\infty}(\|x_n+y_n\|-b)\leq \lim_{n\to\infty}(\|x_n+y_n\|-\|y_n\|)\leq \lim_{n\to\infty}\|x_n\|\leq a
$$
it follows that $\lim_{n\rightarrow \infty}\|x_n\|=a$.

By similar arguments we can prove that $\lim_{n\rightarrow \infty}\|y_n\|=b$.
\end{proof}

\begin{lemma}\label{lem:5}
	Let $(X,\|\cdot\|)$ be a Banach space and let $B=\{x\in X: \|x\|\geq 2\}$. If the ordered pair $(B_X,B)$ satisfies the UC property and the sequences $\{x_n\}_{n=1}^\infty ,\{z_n\}_{n=1}^\infty \subset A$ be such that $\lim_{n\rightarrow \infty}\left\|\frac{x_n+z_n}{2}\right\|= 1$, then 
	$\lim_{n\rightarrow \infty}\|x_n-z_n\|= 0$.
\end{lemma}
\begin{proof}
From the assumptions of the lemma we have that
\begin{equation}\label{l1o1}  
	{\rm dist}(A,B)=1,\ \|x_n\|\leq 1,\ \|z_n\|\leq 1\ \mbox{and}\ \lim_{n\rightarrow \infty}\|x_n+z_n\|= 2.
\end{equation}
Therefore by Lemma \ref{lem:4} it follows
\begin{equation}\label{l1o2}  
	\lim_{n\rightarrow \infty}\|x_n\|= 1,\ \lim_{n\rightarrow \infty}\|z_n\|= 1.
\end{equation}
Let us define the sequence $\{y_n\}_{n=1}^\infty$ by $\displaystyle y_n=2\frac{x_n+z_n}{\|x_n+z_n\|}$. Then $\|y_n\|=2$ for every $n\in\mathbb{N}$, consequently $\{y_n\}_{n=1}^\infty \subset B$. Now using the assumption that $\{x_n\}_{n=1}^\infty\subset A$ we get
the chain of inequalities
\begin{equation}\label{l1o3}  
	\begin{array}{lll} 
	\text{dist}(A,B)&\leq&\|y_n-x_n\|=\left\|2\frac{x_n+z_n}{\|x_n+z_n\|}-x_n\right\|\\[8pt]
	&=&\left\|2\frac{x_n+z_n}{\|x_n+z_n\|}-(x_n+z_n)+z_n\right\|\\[8pt]
		&=&\left\|\left(\frac{2}{\|x_n+z_n\|}-1\right)(x_n+z_n)+z_n\right\|\\
		&\leq& \left(\frac{2}{\|x_n+z_n\|}-1\right)\|x_n+z_n\|+\|z_n\|.
	\end{array}
\end{equation}

From (\ref{l1o1}) and (\ref{l1o2}) it follows that
	$$
	\lim_{n\to\infty}\left(\left(\frac{2}{\|x_n+z_n\|}-1\right)\|x_n+z_n\|+\|z_n\|\right)=1=\text{dist}(A,B)
	$$
	Using (\ref{l1o3}) we get that $\displaystyle\lim_{n\rightarrow \infty}\|y_n-z_n\|={\rm dist}(A,B)$.

By similar arguments we can prove that $\displaystyle\lim_{n\rightarrow \infty}\|y_n-z_n\|={\rm dist}(A,B)$.

From the assumption that the ordered pair $(B_X,B)$ satisfies the $UC$ property it follows that 
$\displaystyle\lim_{n\rightarrow \infty}\|x_n-z_n\|=0$.
\end{proof}

\begin{theorem}\label{Theorem:44}
	Let $(X,\|\cdot\|)$ be a Banach space and let $B=\{x\in X: \|x\|\geq 2\}$. If the ordered pair $(B_X,B)$ satisfies the UC property, then $(X,\|\cdot\|)$ is a uniformly convex Banach space.
\end{theorem}

\begin{proof}
Let us assume {the contrary, i.e.} {there} exist{s} $\varepsilon > 0$ such that for every $\delta > 0$ there are $x(\delta),z(\delta)\in A$, so that the inequalities $\|x(\delta)-z(\delta)\|\geq \varepsilon$ and $1-\delta\leq  \left\|\frac{x(\delta)+z(\delta)}{2}\right\|\leq 1$ hold true. Thus we can choose sequences $\{x_n\}_{n=1}^\infty ,\{z_n\}_{n=1}^\infty \in A$, satisfying $\|x_n-z_n\|\geq \varepsilon$ and
$1-\frac{1}{n}\leq  \left\|\frac{x_n+z_n}{2}\right\|\leq 1$. Therefore $\displaystyle\lim_{n\to \infty}\left\|\frac{x_n+z_n}{2}\right\|= 1$, which contradicts with Lemma \ref{lem:5}.
\end{proof}

\begin{corollary}\label{Corollary:45}
	Let $(X,\|\cdot\|)$ be a Banach space. If every ordered pair of subsets $(A,B)$ has the $UC$ property, where $A$ is a convex, then $(X,\|\cdot\|)$ is a uniformly convex Banach space.
\end{corollary}

\section{The property $UC$ and $UCED$ of the underlying space}

By Proposition \ref{Proposition:18} it follows that in $UCED$ Banach spaces every ordered pair $(A,B)$, such that $A$ is a convex and relatively compact set, satisfies the property $UC$. Unfortunately the sets
$A=\{(x,y)\in (\mathbb{R}^2,\|\cdot\|):y\geq x^2\}$, where $\|\cdot\|$ is a $UCED$ norm or $A=B_X$, where $(X,\|\cdot\|)$ is a $UCED$ Banach space are not relatively compact and therefore for an arbitrary set $B\subset X$ the ordered pair $(A,B)$ does not satisfy the assumption of Proposition \ref{Proposition:18}. The next lemma presents a different condition on the sets of the ordered pair $(A,B)$, by removing the too restrictive assumption the set $A$ to be relatively compact.

\begin{lemma}\label{Lemma:46}
	Let $(X,\|\cdot\|)$ be a Banach space and let $B\subset X/B_X$ such that ${\rm dist}(B,B_X)\geq 1$. If the ordered pair $(B_X,B)$ satisfies the $UC$ property, there exists $p\in B$ so that ${\rm dist}(p,B_X)={\rm dist}(B,B_X)$ and there exist two sequences $\displaystyle\left\{x_n\right\}_{n=1}^\infty ,\left\{z_n\right\}_{n=1}^\infty \subset B_X$, satisfying $\displaystyle\lim_{n\to\infty}\left\|\frac{x_n+z_n}{2}\right\|= 1$ and $x_n+z_n=|\lambda_n| p$, then 
	$\lim_{n\rightarrow \infty}\|x_n-z_n\|= 0$.
\end{lemma}
\begin{proof}
From the assumptions we have that
$\|x_n\|\leq 1$ , $\|z_n\|\leq 1$ and $\lim_{n\rightarrow \infty}\|x_n+z_n\|= 2$. Therefore from Lemma \ref{lem:4} it follows that 
\begin{equation}\label{equation:27}  
	\lim_{n\rightarrow \infty}\|x_n\|= 1\ ,\ \lim_{n\rightarrow \infty}\|z_n\|= 1\ \mbox{and}\ \lim_{n\rightarrow \infty}\|x_n+z_n\|= 2.
\end{equation}
Let us define the sequence $\{y_n\}_{n=1}^\infty$ by $y_n=p$ for every $n\in\mathbb{N}$.

From the assumption $x_n+z_n=|\lambda_n| p$ it follows that $\|p-(x_n+z_n)\|=|\|p\|-\|x_n+z_n\||$. 

After using the inequalities $\|x_n+z_n\|\leq \|x_n\|+\|z_n\|\leq 2$, ${\rm dist}(B,B_X)+1\geq 2$, the equality
$\|p\|={\rm dist}(p,B_X)+1={\rm dist}(B,B_X)+1$ we can write the chain of inequalities
\begin{equation}\label{eeeee}  
	\begin{array}{lll}
		{\rm dist}(B,B_X)&\leq&\displaystyle\|y_n-x_n\|=\|p-x_n\|=\|p-(x_n+z_n)+z_n\|\\
		&\leq& \displaystyle\left(\|p-(x_n+z_n)\|+\|z_n\|\right)\\
		&=&\displaystyle\left(|\|p\|-\|x_n+z_n\|| +\|z_n\|\right)\\
		&=&\displaystyle\left(|{\rm dist}(B,B_X)+1-\|x_n+z_n\|| +\|z_n\|\right)\\
		&=& {\rm dist}(B,B_X)+1+\displaystyle\left(\|z_n\|-\|x_n+z_n\|\right).
	\end{array}
\end{equation}

From (\ref{equation:27}) we get 
	$$\lim_{n\to\infty}\left({\rm dist}(B,B_X)+1+\displaystyle\|z_n\|-\|x_n+z_n\|\right)={\rm dist}(B,B_X).$$
	Thus using (\ref{eeeee}) it follows that
Therefore $\displaystyle\lim_{n\to\infty}\|y_n-x_n\|={\rm dist}(B_X,B)$. 

By similar arguments we can prove that 
$\displaystyle\lim_{n\to\infty}\|y_n-z_n\|={\rm dist}(B_X,B)$.

From the assumption that the ordered pair $(B_X,B)$ satisfies the $UC$ property it follows that 
$\displaystyle\lim_{n\to\infty}\|x_n-z_n\|=0$.
\end{proof}

{\begin{definition}\label{Definition:47}
		We say that a Banach space $(X,\|\cdot\|)$ is uniformly convex in the direction $z$ if
		$\delta_{\|\cdot\|}(z,\varepsilon)>0$ for every $\varepsilon\in (0,2]$.
\end{definition}}

\begin{theorem}\label{Theorem:48}
	Let $(X,\|\cdot\|)$ be a Banach space and let $B\subset X/B_X$. If there holds ${\rm dist}(B,B_X)\geq 1$, the ordered pair $(B_X,B)$ satisfies the $UC$ property and there is $p\in B$ so that ${\rm dist}(p,B_X)={\rm dist}(B,B_X)$, then $(X,\|\cdot\|)$ is a uniformly convex {in} the direction $\displaystyle\frac{p}{\|p\|}$.
\end{theorem}

\begin{proof}
Let {us} assume {the contrary, i.e.} {there} exists $\varepsilon > 0$ such that for every $\delta > 0$ there are $x(\delta),z(\delta)\in B_X$, satisfying $\|x(\delta)-z(\delta)\|\geq \varepsilon$, $\displaystyle 1-\delta\leq  \left\|\frac{x(\delta)+z(\delta)}{2}\right\|\leq 1$ and  $x(\delta)+z(\delta)=\lambda(\delta) p$. Thus we can choose sequences $\{x_n\}_{n=1}^\infty ,\{z_n\}_{n=1}^\infty \in B_X$ so that $\|x_n-z_n\|\geq \varepsilon$, $\displaystyle\lim_{n\to \infty}\left\|\frac{x_n+z_n}{2}\right\|= 1$ and $x_n+z_n=\lambda_n p$. 

Let us construct the sequences $\{x'_n\}_{n=1}^\infty$,$\{z'_n\}_{n=1}^\infty\subset B_X$ as follows: $x'_n=x_n$, $z'_n=z_n$ if $\lambda_n\geq 0$ and $x'_n=-x_n$, $z'_n=-z_n$ if $\lambda_n<0$. Then $\{x'_n\}_{n=1}^\infty\,\{z'_n\}_{n=1}^\infty\in B_X$, $\|x'_n-z'_n\|\geq \varepsilon$, $\displaystyle\lim_{n\to\infty}\left\|\frac{x'_n+z'_n}{2}\right\|= 1$ and $x'_n+z'_n=|\lambda_n| p$ which contradicts with Lemma \ref{Lemma:46} and consequently  $(X,\|\cdot\|)$ is a uniformly convex in the direction $\displaystyle\frac{p}{\|p\|}$. 
\end{proof}

\section{Examples and applications}

We will finish with an example of a cyclic map, which satisfies the conditions of Theorem \ref{Theorem:14}, but the ordered pair of sets
$(A,B)$, that are the domain of the map $T$, are just $BUC$ but not $UC$. Therefore we can not conclude the existence of best proximity points using the result of \cite{Suzuki-Kikkawa-Vetro}.

\begin{example}\label{Example:49}
	Let us consider the space $X=(\mathbb{R}^2,\|\cdot\|_\infty)$, the sets 
	$$\displaystyle A=\left\{(x,y): \frac{1}{x}\leq y, x>0 \right\}$$
	and
	$$
	B=B_1\cup B_2,
	$$
	where
	$$
	\displaystyle B_1=\left\{(x,y):\frac{1}{x+1}-1\geq y\ and \ x>-1\right\}
	$$ 
	and
	$$\displaystyle B_2=\left\{(x,y):y\in\mathbb{R}, x\leq -1\right\}.
	$$
	Let us denote 
	$$
	\displaystyle\overline{A}=\left\{(x,y): \frac{1}{x}=y, x>0 \right\},
	$$
	and
	$$
	\displaystyle \overline{B}=\left\{(x,y):\frac{1}{x+1}-1= y\ and \ x>-1\right\},
	$$
	Let the map $T:A\cup B\rightarrow A\cup B$ be defined by 
	\begin{equation}\label{equation:28}
		Tx=\left\{
		\begin{array}{rl}
			\displaystyle\left(-\frac{{\rm dist}(x,\overline{A})}{2},-\frac{{\rm dist}(x,\overline{A})}{2}\right),& x\in A\\[16pt]
			\displaystyle\left(1+\frac{{\rm dist}(x,\overline{B})}{2},1+\frac{{\rm dist}(x,\overline{B})}{2}\right),& x\in B.
		\end{array}
		\right.
	\end{equation} 
	Then $T$ is a cyclic map on $A\cup B$ and satisfies the inequality
	\begin{equation}\label{equation:29}
		\rho(Tx,Ty)\leq\frac{1}{2}\rho(x,y)+\frac{1}{2}{\rm dist}(A,B)
	\end{equation} 
	for every $x\in A$ , $y\in B$, i.e. Theorems \ref{Theorem:31} with $\displaystyle k=\frac{1}{2}$.
\end{example}

We will prove first that $T$ is a cyclic map, i.e. $T(A)\subseteq B$ and $T(B)\subseteq A$.

Let first $z=(x,y)\in A$. Then $\displaystyle Tz=\displaystyle\left(-\frac{{\rm dist}(z,\overline{A})}{2},-\frac{{\rm {\rm dist}}(z,\overline{A})}{2}\right)$. 

By $\displaystyle -\frac{{\rm dist}(z,\overline{A})}{2}\leq 0$ it follows that for any 
$z\in A$, if $\displaystyle-\frac{{\rm dist}(z,\overline{A})}{2}\leq -1$ there holds the inclusion $\displaystyle\left(-\frac{{\rm dist}(z,\overline{A})}{2},-\frac{{\rm dist}(z,\overline{A})}{2}\right)\in B_2$. If
$\displaystyle-\frac{{\rm dist}(z,\overline{A})}{2}\in (-1,0]$, then $\displaystyle\left(-\displaystyle\frac{{\rm dist}(z,\overline{A})}{2},-\frac{{\rm dist}(z,\overline{A})}{2}\right)\in B_1$. Thus $Tz\in B$ for every $z\in A$.

Let $z=(x,y)\in B$. Then $\displaystyle Tz=\left(1+\frac{{\rm dist}(z,\overline{B})}{2},1+\frac{{\rm dist}(z,\overline{B})}{2}\right)$. 

From $\displaystyle 1+\frac{{\rm dist}(z,\overline{B})}{2}\geq 1$ and $\displaystyle \frac{1}{x}\leq 1$ for $x\geq 1$ it follows 
$$
\displaystyle\left(1+\frac{{\rm dist}(z,\overline{B})}{2},1+\frac{{\rm dist}(z,\overline{B})}{2}\right)\in A,
$$ 
i.e. $Tz\in A$.  

It remains to show that the map $T$ satisfies the inequality (\ref{equation:29}).

It is easy to calculate that
$$
{\rm dist}(A,B)=\inf\{\|x-y\|_\infty:x\in A,y\in B\}=1
$$
and
$$
{\rm dist}(a,\overline{A})+{\rm dist}(A,B)+{\rm dist}(b,\overline{B})\leq \|a-b\|_\infty.
$$

Let $a=(x_a,y_a)\in A$ and $b=(x_b,y_b)\in B$, then  
$$
\begin{array}{lll}
	\textstyle\|Ta-Tb\|_\infty&=&\displaystyle\left\|\left(1+\textstyle\frac{{\rm dist}(b,\overline{B})}{2},1+\frac{{\rm dist}(b,\overline{B})}{2}\right)-\left(-\textstyle\frac{{\rm dist}(a,\overline{A})}{2},-\frac{{\rm dist}(a,\overline{A})}{2}\right)\right\|_\infty\\
	&=&\displaystyle1+\frac{{\rm dist}(b,\overline{B})}{2}+\frac{{\rm dist}(a,\overline{A})}{2}\\
	&=&\displaystyle\frac{1}{2}\left({\rm dist}(A,B)+{\rm dist}(b,\overline{B})+{\rm dist}(a,\overline{A})\right)+\frac{1}{2}{\rm dist}(A,B)\\[8pt]
	&\leq&\displaystyle \frac{1}{2}\|a-b\|_\infty+\frac{1}{2}{\rm dist}(A,B).
\end{array}
$$ 

The the ordered pair of sets $(A,B)$ has the $BUC$ property, but not the $UC$ property by Example \ref{Example:43} therefore by Theorem \ref{Theorem:33} it follows that there is a unique best proximity point of $T$ in $A$. As far as $(A,B)$ has not the $UC$ property we can not use Theorem \ref{Theorem:14} to conclude that there is a best proximity point of $T$ in $A$.

We have seen that in order the key lemma from \cite{Eldred-Veeramani} to hold true we need either to impose good properties on the unit ball $B_X$ of the underlying Banach space $(X,\|\cdot\|)$ or to impose good properties on the set $A$, from the ordered pair $(A,B)$. The set $A$ in Example \ref{Example:49} is a strictly convex set, i.e. for every $x,y\in A$, there is $\varepsilon=\varepsilon(x,y)>0$ so that $B\left(\frac{x+y}{2},\varepsilon\right)\subset A$ and we have proven that any ordered pair $(A,B)$ has the $BUC$ property. 

We will show in the next example, where the underlying space is strictly convex without being uniformly convex. Thereafter we can not apply Lemma \ref{Lemma:6}. We will show that the ordered pair $(B_X,\{x\in X:\|x\|\geq 2\})$ has not either $BUC$ nor $UC$.

\begin{example}\label{Example:50}
	Let us denote $E_n=(R^2,\|\cdot\|_n)$, where $\|(x,y)\|_n=\root {n} \of {|x|^n+|y|^n}$.
	Any of the spaces $E_n$ is a uniformly convex Banach space. Let us consider the space $\displaystyle X=\left(\prod_{n=2}^\infty E_n,\|\cdot\|\right)$, where $\displaystyle\|\{x_n\}\|=\sqrt{\sum_{n=2}^\infty \|x_n\|_n^2}$ and $x_n\in E_n$. The space $X$ is a strictly convex Banach space, which is not uniformly convex, i.e. its unit ball $B_X$ is strictly convex set. 
	
	Let us denote $B=\{x\in X:\|x\|\geq 2\}$. We will construct three sequences $\left\{x_n\right\}_{n=1}^{\infty},\{z_n\}_{n=1}^{\infty}\in B_x$ and $\left\{y_n\right\}_{n=1}^{\infty}\in B$ such that $\displaystyle\lim_{n\rightarrow \infty}\|x_n-y_n\|= {\rm dist}(B_X,B)$, $\displaystyle\lim_{n\rightarrow \infty}\|z_n-y_n\|= {\rm dist}(B_X,B)$ and $\lim_{n\rightarrow \infty}\|x_n-z_n\|> 0$, i.e. the ordered pair $(B_X,B)$ has not the $UC$ property.
\end{example}

Let $r_x,r_z:\mathbb{N}\rightarrow \mathbb{R}^2$ be defined by $r_x(n)=\left(\frac{1}{\sqrt[n]{2}},\frac{1}{\sqrt[n]{2}}\right)$, $r_z(n)=\left(\frac{1}{\sqrt[n]{2}},-\frac{1}{\sqrt[n]{2}}\right)$. 

We can see that $\|r_x(n)\|_n = \|r_z(n)\|_n = 1$. 

Let 
$$
x_n=(\underbrace{(0,0),(0,0),...,(0,0)}_{n-1},r_x(n+1),(0,0),...),
$$
$$z_n=(\underbrace{(0,0),(0,0),...,(0,0)}_{n-1},r_z(n+1),(0,0),...)
$$ 
and 
$$y_n=(\underbrace{(0,0),(0,0),...,(0,0)}_{n-1},(2,0),(0,0),...).
$$ 
From $\|x_n\|=\|r_x(n+1)\|_{n+1}=1$, $\|z_n\|=\|r_z(n+1)\|_{n+1}=1$ and $\|y_n\|=2$ it follows that $\left\{x_n\right\}_{n=1}^{\infty},\left\{z_n\right\}_{n=1}^{\infty}\subset B_X$ and $\left\{y_n\right\}_{n=1}^{\infty}\subset B$

For every $n\in\mathbb{N}$ there holds
$$
\begin{array}{lll}
	\|x_n-y_n\|&=&\displaystyle\|r_x(n+1)-(2,0)\|_{n+1}=\left\|\left(\frac{1}{\sqrt[n+1]{2}}-2,\frac{1}{\sqrt[n+1]{2}}\right)\right\|_{n+1}\\
	&=&\displaystyle\sqrt[n+1]{\left(2-\frac{1}{\sqrt[n+1]{2}}\right)^{n+1}+\frac{1}{2}}.
\end{array}
$$
Using the inequalities
$$
\begin{array}{lll}
	2-\frac{1}{\sqrt[n+1]{2}}&=&\displaystyle\sqrt[n+1]{\left(2-\frac{1}{\sqrt[n+1]{2}}\right)^{n+1}}\leq\sqrt[n+1]{\left(2-\frac{1}{\sqrt[n+1]{2}}\right)^{n+1}+\frac{1}{2}}\\
	&<&\displaystyle\sqrt[n+1]{2\left(2-\frac{1}{\sqrt[n+1]{2}}\right)^{n+1}}=2\sqrt[n+1]{2}-1
\end{array}
$$
for every $n\in\mathbb{N}$, we get $\displaystyle\lim_{n\to\infty}\|x_n-y_n\|=1$. 

We can prove that $\displaystyle\lim_{n\to\infty}\|z_n-y_n\|=1$ in a similar fashion.

From $\displaystyle\|x_n-z_n\|=\|r_x(n+1)-r_z(n+1)\|_{n+1}=\left\|\left(0,\frac{2}{\sqrt[n+1]{2}}\right)\right\|_{n+1}=\frac{2}{\sqrt[n+1]{2}}$ for every $n\in\mathbb{N}$ it follows that $\displaystyle\lim_{n\rightarrow \infty}\|x_n-z_n\|=2$.

Thus $(B_X,B)$ has not the $UC$ property.
The sequences in the last example are bounded ones, and therefore ordered pair $(B_X,B)$ has not the $BUC$ property, too. 

In view of Theorem \ref{Theorem:42}, $B_X$ is not uniformly convex about any function $\phi$. As far as $B_X$ is a strictly convex set it follows that there is a difference between a strict convexity of a set and a uniform convexity of a set about a function $\varphi$.

\end{document}